\documentclass{article}
\usepackage{amsmath}
\usepackage{amssymb}

\newcommand{\T}{\mathbb{T}}
\newcommand{\D}{\mathbb{D}}
\newcommand{\C}{\mathbb{C}}
\newcommand{\R}{\mathbb{R}}
\newcommand{\Z}{\mathbb{Z}}
\newcommand{\N}{\mathbb{N}}

\def\co{\colon}
\def\ba{^{\relbar}}
\def\la{\langle}
\def\ra{\rangle}
\def\i{\infty}

\def\op{\oplus}
\def\eq{\simeq}
\def\ol{\overline}

\def\HH{{\cal H}}
\def\L{{\cal L}}
\def\M{{\cal M}}
\def\K{{\cal K}}
\def\I{{\cal I}}
\def\TT{{\cal T}}
\def\P{{\cal P}}
\def\A{{\cal A}}
\def\X{{\cal X}}
\def\B{{\cal B}}

\def\a{\alpha}
\def\b{\beta}
\def\g{\gamma}
\def\f{\varphi}
\def\ff{\phi}
\def\o{\omega}
\def\O{\Omega}

\def\z{\zeta}

\def\s{\sigma}

\def\FF{\Phi}

\def\h{\chi}

\def\m{\mu}

\def\ae{\overset{a}\sim}
\def\wt{\widetilde}

\newtheorem{pro}{Proposition}
\newtheorem{thm}[pro]{Theorem}
\newtheorem{lem}[pro]{Lemma}
\newtheorem{cor}[pro]{Corollary}

\newtheorem{rem}[pro]{Remark}
\newtheorem{que}[pro]{Question}

\begin{document}

\title{Quasianalytic polynomially bounded operators}

\author{L\'aszl\'o K\'erchy \\
University of Szeged}

\date{}

\maketitle

\begin{abstract}
\noindent
Quasianalytic contractions form the crucial class in the quest for proper invariant and hyperinvariant subspaces for asymptotically non-vanishing Hilbert space contractions.
The property of quasianalycity relies on the concepts of unitary asymptote and $H^\i$-functional calculus.
These objects can be naturally defined in the setting of polynomially bounded operators too, which makes possible to extend the study of quasianalycity from contractions to this larger class.
Carrying out this program we pose also several interesting questions.


\end{abstract}

\section{Introduction} 
\label{intr}

Let $\HH$ be a complex Hilbert space, and let $\L(\HH)$ stand for the $C^*$-algebra of all bounded linear operators acting on $\HH$.
Given any $T\in\L(\HH)$, the complete lattice $\hbox{ Lat} \,T$ consists of those (closed) subspaces $\M$ of $\HH$, which are \emph{invariant} for $T: \; T\M\subset\M$.
A subspace $\M$ is \emph{hyperinvariant} for $T$, if it is invariant for every operator $C\in\L(\HH)$, commuting with $T: \; CT=TC$.
The complete lattice of all hyperinvariant subspaces of $T$ is denoted by $\hbox{ Hlat}\, T$.
The most challenging open questions in operator theory are arguably the \emph{Invariant Subspace Problem} (ISP) and the \emph{Hyperinvariant Subspace Problem} (HSP).
The first question asks whether $\hbox{ Lat}\,T$ is non-trivial (i.e., different from $\{\{0\},\HH\}$), for every operator $T\in\L(\HH)$, provided dim $\HH \ge 2$; while the second question asks whether $\hbox{ Hlat}\, T$ is non-trivial, whenever $T$ is not a scalar multiple of the identity operator $I$.
For a thorough discussion of these problems see the monographs \cite{Radjavi} and \cite{Challendar}, and the more recent papers \cite{Argyros} and \cite{Grivaux}.

Of course, dividing $T$ by its norm, we may assume that $T$ is a contraction: $\|T\|\le 1$.
An opposite, local estimate is made, when we assume that the contraction $T$ is \emph{asymptotically non-vanishing}, that is $\lim_n\|T^nh\|>0$ holds for some $h\in\HH$.
Suprisingly enough, (ISP) and (HSP) are open under these assumptions too.
A useful extra tool available in this situation is the unitary asymptote of $T$.
Relying on this and on the $H^\i$-functional calculus, two spectral invariants can be introduced for $T$.
The first one is the residual set $\o(T)$, and the second one is the quasianalytic spectral set $\pi(T)$ of $T$.
The contraction $T$ is called quasianalytic, if these measurable subsets of the unit circle $\T$ coincide: $\o(T)=\pi(T)$.
Quasianalytic contractions were investigated in the papers \cite{Ker2001}, \cite{Ker2011}, \cite{Ker2013}, \cite{KerTot}, \cite{KSzStudia} and \cite{KSzProc}.
A central theorem proved is that (HSP) can be reduced to this class in the asymptotically non-vanishing case.

The aim of this paper is to extend these investigations from contractions to polynomially bounded operators.
It will turn out that this larger class is the natural setting for the study of quasianalycity.

Our work is organized in the following way.

In Section \ref{unitasymp} we introduce unitary asymptotes for an arbitrary operator $T\in\L(\HH)$ in a categorical sense, as it has been done for contractions in \cite{BK} and in Chapter IX of \cite{NFBK}.
The induced generalized Toeplitz operators, and the connection of the associated symbolic calculus with the commutant mapping are discussed.
Existence of unitary asymptotes in the power bounded case is shown by using Banach limits.
This method originates in Sz.-Nagy's pioneering work \cite{Nagy47}, and has been extended to many situations; see, e.g., \cite{KerInd89}, \cite{KerNormseq}, \cite{KerDiscrepr}, \cite{KerProc}, \cite{KerToepl},  \cite{KerLeka},  \cite{KerMul} and \cite{Prunaru}.
A new characterization is given in terms of norm-conditions.

It is well-known that every contraction $T$ can be decomposed into the orthogonal sum $T= T^a \op U^s$ of an absolutely continuous contraction $T^a$ and a singular unitary operator $U^s$.
Mlak showed that analogous Lebesgue-type decompositions can be given for polynomially bounded operators, or equivalently, for bounded representations of the disk algebra $\A(\T)$. 
Actually, Mlak considered representations of more general function algebras; see \cite{Mlak Acta}, \cite{Mlak singular}, \cite{Mlak Oberwolfach} and \cite{Mlak algebraic}.
Using Mlak's elementary measures, in Section \ref{decomp} we give a detailed, streamlined discussion of the Lebesgue decomposition in the particular case $\A(\T)$, what is of the main interest for us.
Our purpose here (and partly in the next two sections) is to make Mlak's important results more accessible for a wide range of readers.

In Section \ref{intertw} we focus on intertwining relations.
It is verified that absolute continuity and singularity are preserved under quasisimilarity.
A transparent proof is given for the important known fact that every singular polynomially bounded operator is similar to its unitary asymptote.
Furthermore, the problem of similarity to contractions is discussed.

Sz.-Nagy and Foias introduced the effective functional calculus, working with functions in $H^\i$, for absolutely continuous contractions; see Chapter III in \cite{NFBK}.
Section \ref{calculus} is devoted to the study of properties and possible range of this $H^\i$-functional calculus.
A simplified proof is given for Mlak's result, stating that exactly the absolutely continuous polynomially bounded operators admit $H^\i$-functional calculae.

In Section \ref{quasianalytic} we turn to quasianalytic operators.
The sets $\o(T)$ and $\pi(T)$ are introduced in a uniform way, relying on the local residual sets $\o(T,x)\; (x\in\HH)$.
The quasianalytic spectral set $\pi(T)$ is characterized also in terms of non-monotone sequences.
The central hyperinvariant subspace theorem is proved.
Furthermore, the effects of intertwining relations, the asymptotic behaviour, orthogonal sums and restrictions are studied.

For the theory of Hardy spaces $H^p$, see \cite{Hof}.
In connection with the theory of contractions we refer to \cite{NFBK}.
Given operators $T_1\in\L(\HH_1)$ and $T_2\in\L(\HH_2)$, the set of \emph{intertwining transformations} is defined by
\[\I(T_1,T_2) := \left\{ C\in\L(\HH_1,\HH_2): CT_1 = T_2C\right\}.\]
The \emph{commutant} of the operator $T\in\L(\HH)$ is $\{T\}'= \I(T,T)$, while the \emph{bicommutant} (or double commutant) $\{T\}''$ of $T$ consists of those operators $B\in\L(\HH)$ which commute with every operator $C$ in $\{T\}'$.
Finally, $\N, \Z, \Z_+, \R, \R_+, \C$ denote the set of positive integers, integers, non-negative integers, real numbers, non-negative real numbers and complex numbers, respectively.

\section{Unitary asymptotes}
\label{unitasymp}

Let $\HH$ be a complex Hilbert space, and let us consider an operator $T\in\L(\HH)$.
We say that $(X,U)$ is a \emph{unitary intertwining pair} for $T$, if $U$ is a unitary operator acting on a (complex) Hilbert space $\K$, and $X\in\I(T,U)$.
The pair is \emph{minimal}, if $\vee_{n=1}^\i U^{-n}X\HH= \K$.
The unitary intertwining pair $(X,U)$ is called a \emph{unitary asymptote} of $T$, if for any other unitary intertwining pair $(X',U')$ of $T$, there exists a unique transformation $Z\in\I(U,U')$ such that $X'=ZX$.
The uniqueness of $Z$ implies that a unitary asymptote $(X,U)$ is necessarily minimal.

Let us assume that $(X_1,U_1)$ and $(X_2,U_2)$ are unitary asymptotes of $T$.
Then there exist $Z_1\in\I(U_1,U_2)$ and $Z_2\in\I(U_2,U_1)$ such that $Z_1X_1=X_2$ and $Z_2X_2=X_1$.
Therefore, $I\cdot X_1=Z_2X_2= (Z_2Z_1)X_1$, whence $I=Z_2Z_1$ follows by the definition of the unitary asymptote.
The equation $Z_1Z_2=I$ can be derived similarly.
Thus, the unitary intertwining pairs $(X_1,U_1)$ and $(X_2,U_2)$ are \emph{equivalent}, which means the existence of an invertible transformation $Z\in\I(U_1,U_2)$ satisfying $X_2=ZX_1$.
Then the unitary operators $U_1$ and $U_2$ are necessarily unitarily equivalent: $U_1\eq U_2$.
Furthermore, if the unitary intertwining pair $(X_2,U_2)$ is equivalent to a unitary asymptote $(X_1,U_1)$, then $(X_2,U_2)$ is also a unitary asymptote.

For any unitary intertwining pair $(X,U)$ of $T$ we have
$ \|Xh\|= \|U^nXh\| = \|XT^nh\|\le \|X\|\cdot \|T^nh\| \; (n\in\N)$, and so
\[ \|Xh\|\le \|X\|\cdot \liminf_{n\to\i}\|T^nh\|\quad \hbox{ for all }\; h\in\HH.\]
A lower estimate for $\|Xh\|$ yields universality.

\begin{pro}
\label{suff}
Let $(X,U)$ be a minimal unitary intertwining pair for $T$.
If there exists $\kappa \in (0,\i)$ such that
\[ \kappa\cdot \liminf_{n\to\i} \|T^nh\| \le \|Xh\| \]
holds for all $h\in\HH$, then $(X,U)$ is a unitary asymptote of $T$.
\end{pro}

\noindent
\textbf{Proof.}
Let $(X',U')$ be a unitary intertwining pair for $T$.
For every $h\in\HH$ we have
\[ \|X'h\| \le \|X'\|\cdot \liminf_{n\to\i} \|T^nh\| \le \|X'\|\cdot \kappa^{-1}\cdot \|Xh\|.\]
Hence there is a unique transformation $Y_+\in \L(\K_+,\K')$ such that $Y_+Xh=X'h\; (h\in\HH)$, where $\K_+= (X\HH)\ba$.
The relations
\[ Y_+UXh=Y_+XTh=X'Th=U'X'h=U'Y_+Xh\quad (h\in \HH)\]
show that $Y_+\in \I(U_+,U')$, where $U_+=U|\K_+$.
It is immediate that the equations
\[YU^{-n}k:= {U'}^{-n}Y_+k \quad (k\in \K_+, n\in\N)\]
define a norm-preserving extension $Y\in \I(U,U')$ of $Y_+$.
It is clear that $YX=X'$, and that $Y$ is uniquely determined by these properties.

\rightline{$\square$}

\begin{que}
\label{necessary}
Does $\kappa\liminf_{n\to\i}\|T^nh\|\le \|Xh\|\; (h\in\HH)$ hold with a $\kappa\in(0,\i)$, whenever $(X,U)$ is a unitary asymptote of $T$?
\end{que}

Let $(X,U)$ be any unitary intertwining pair for $T$.
Then $T^*(X^*X)T= X^*U^*UX= X^*X$.
Thus $X^*X$ belongs to the set 
\[ \TT(T):= \{B\in \L(\HH) : T^*BT=B\}\]
of \emph{$T$-Toeplitz operators}.
It is clear that $\TT(T)$ is a selfadjoint linear manifold, which is closed in the weak operator topology.
Taking $\TT_s(T)= \{B\in\TT(T) : B^*=B\}$, we have $\TT(T)= \TT_s(T)\; \dot+ \;i\TT_s(T)$. 
The set $\TT_+(T)= \{B\in\TT_s(T) : B\ge 0\}$ is of particular interest.
We recall that $A\in \TT_+(T)$ is \emph{universal} (\emph{weakly universal}) in $\TT(T)$, if for every $B\in \TT_s(T)$ ($B\in \TT_+(T)$, respectively) there exists $\b\in\R_+$ such that
\[-\b A\le B \le \b A.\]
It is obvious that $\TT(U)=\{U\}'$.
Furthermore, for any $F\in \{U\}'$ we have $ T^*(X^*FX)T = X^*U^*FUX = X^*U^*UFX = X^*FX$, that is $X^*FX\in\TT(T)$.
The positive, bounded, $*$-linear mapping
\[\Psi_X\co \{U\}' \to \TT(T), \; F\mapsto X^*FX\]
is the \emph{symbolic calculus} for $\TT(T)$, associated with $(X,U)$.
It is easy to see that $\Psi_X$ is injective, when $(X,U)$ is minimal.
We say that $\Psi_X$ is \emph{universal}, if its range coincides with $\TT(T)$.
The following theorem, proved in \cite{KerToepl} in a slightly different form, establishes connections among these concepts.

\begin{thm}
\label{connections}
Let $(X,U)$ be a minimal unitary intertwining pair for $T$.
Then $\textup{(i)} \iff \textup{(ii)} \Leftarrow \textup{(iii)} \iff \textup{(iv)}$, where
\begin{itemize}
\item[\textup{(i)}] $(X,U)$ is a unitary asymptote of $T$;
\item[\textup{(ii)}] $A=X^*X\in \TT_+(T)$ is weakly universal in $\TT(T)$;
\item[\textup{(iii)}] $A=X^*X\in\TT_+(T)$ is universal in $\TT(T)$;
\item[\textup{(iv)}] the symbolic calculus $\Psi_X$ is universal.
\end{itemize}
\end{thm}

We note that if $A=X^*X$ is universal in $\TT(T)$, then $\Psi_X\co \{U\}' \to \TT(T)$ is invertible, and so
\[ - \b_0\|B\| A \le B\le \b_0\|B\| A\]
holds for every $B\in \TT_s(T)$, where $\b_0=\|\Psi_X^{-1}\|.$

It may happen that there is not a universal operator in $\TT(T)$; we refer to Example 5 in \cite{KerToepl}.
Moreover, if $\TT(T)$ contains a non-zero operator $B$, then the relations $0< \|B\|= \|T^{*n}BT^n\|\le \|B\| \|T^n\|^2\; (n\in\N)$ show that the spectral radius $r(T)\ge 1$.

Let us assume that $(X,U)$ is a unitary asymptote of $T$.
Then, for every $C\in\{T\}', XC\in\I(T,U)$ holds, and so there exists a unique $D\in\{U\}'$ such that $XC=DX$.
It is easy to verify that the \emph{commutant mapping}
\[\g_X\co \{T\}' \to \{U\}', \; C\mapsto D\]
is a (unital) algebra-homomorphism.
Boundedness of $\g_X$ follows from its relation to $\Psi_X$.

\begin{pro}
\label{commutant mapping}
Let $(X,U)$ be a unitary asymptote of $\,T$, and let us assume that $A=X^*X$ is universal in $\TT(T)$.
Then $(\Psi_X\circ\g_X)C= AC$ holds for every $C\in\{T\}'$, and so $\|\g_X\|\le \|A\| \cdot \|\Psi_X^{-1}\|.$
\end{pro}

The simple proof is left to the reader.

\medskip

Existence of unitary asymptotes can be verified in the class of power bounded operators.
We recall that $T\in\L(\HH)$ is a \emph{power bounded operator}, if
\[ M_T:= \sup\{ \|T^n\| : n\in \Z_+\}< \i.\]
In that case $r(T)\le1$ is evidently true.
If $r(T)<1$, then $\lim_{n\to\i} \|T^n\|=0$, and so, for any unitary intertwining pair $(X,U)$ of $T$, the equations $XT^n=U^nX\; (n\in\N)$ impõly that $X=0$.
Thus $U$ acts on the zero space, provided $(X,U)$ is minimal.
Therefore $r(T)=1$, whenever the power bounded operator $T$ has a non-trivial unitary asymptote.

Taking a Banach limit $L\co  \ell^\i \to\C$, there exists a unique operator $A_L\in\L(\HH)$ such that
\[ \la A_Lx,y\ra = L\hbox{-}\lim_{n\to\i} \la T^{*n}T^nx,y\ra \quad \hbox{ for all } \; x,y\in\HH.\]
It is evident that $A_L\in \TT_+(T)$, and $A_L\le M_T^2\cdot I.$ 
Setting any $B\in \TT_s(T)$, the relations $\la Bx,x\ra = \la BT^nx,T^nx\ra$ and
\[ - \|B\| \la T^nx,T^nx\ra \le \la BT^nx,T^nx\ra \le \|B\| \la T^nx,T^nx\ra \quad (x\in\HH, n\in\N)\]
show that
\[ -\|B\| A_L \le B \le \|B\| A_L.\]
Therefore, $A_L$ is universal in $\TT(T)$.

Since $T^*A_LT=A_L$ implies $\|A_L^{1/2}Th\|= \|A_L^{1/2}h\|\; (h\in\HH)$, there exists a unique isometry $V_L$ on the space $\K_{+,L}= (A_L^{1/2}\HH)\ba$ satisying the condition $V_LA_L^{1/2}=A_L^{1/2}T$.
Let $U_L\in \L(\K_L)$ be the minimal unitary extension of $V_L$, and let $X_L\in \L(\HH,\K_L)$ be defined by $X_Lh:= A_L^{1/2}h\; (h\in\HH)$.
Clearly, $(X_L,U_L)$ is a minimal unitary intertwining pair for $T$.
Taking into account that $X_L^*X_L= A_L$ is universal in $\TT(T)$, we infer by Theorem \ref{connections} that $(X_L,U_L)$ \emph{is a unitary asymptote of $T$}.

We recall that the possible values of Banach limits on a real sequence $\xi\in\ell^\i$ comprise a closed interval determined by the numbers
\[\underline B(\xi)= \underline B\hbox{-}\lim_{n\to\i}\xi(n) := \sup\left\{ \liminf_{k\to\i} \frac{1}{r}\sum_{j=1}^r \xi(n_j+k) : n_1,\dots,n_r\in\N, r\in\N\right\} \]
and
\[\ol B(\xi)= \ol B\hbox{-}\lim_{n\to\i}\xi(n) := \inf\left\{ \limsup_{k\to\i} \frac{1}{r}\sum_{j=1}^r \xi(n_j+k) : n_1,\dots,n_r\in\N, r\in\N\right\}.\] 
For any $h\in\HH$, we have
\[\|X_Lh\|^2=\la A_Lh,h\ra = L\hbox{-}\lim_{n\to\i}\|T^nh\|^2,\]
whence
\[\underline B\hbox{-}\lim_{n\to\i}\|T^nh\|^2\le \|X_Lh\|^2\le \ol B\hbox{-}\lim_{n\to\i}\|T^nh\|^2.\]
These relations immediately yield that
\[\liminf_{n\to\i}\|T^nh\|\le \|X_Lh\|\le \limsup_{n\to\i}\|T^nh\|\quad (h\in\HH).\]
We obtain the following characterization.

\begin{thm}
\label{characterization}
Let $(X,U)$ be a minimal unitary intertwining pair for the power bounded operator $T\in\L(\HH)$.
Then $(X,U)$ is a unitary asymptote of $T$ if and only if there exists a $\kappa\in (0,\i)$ such that
\[\kappa\cdot \liminf_{n\to\i} \|T^nh\| \le \|Xh\| \quad \hbox{ for all } \; h\in\HH.\]
\end{thm}

\noindent
\textbf{Proof.}
If $(X,U)$ is a unitary asymptote of $T$, then there exists an invertible $Z\in\I(U,U_L)$ such that $X_L=ZX$, with a Banach limit $L$.
Thus, for every $h\in\HH$, we have
\[\| Xh\| \ge \|Z\|^{-1} \|ZXh\|= \|Z\|^{-1} \|X_Lh\|\ge \|Z\|^{-1}\liminf_{n\to\i}\|T^n h\|.\]
The reverse implication follows by Proposition \ref{suff}.

\rightline{$\square$}

As immediate consequences, we may derive the following statements.

\begin{cor}
\label{direct sum}
Let $(X_j,U_j)$ be a unitary asymptote of the power bounded operator $T_j$, for $j=1,2.$
Then $(X_1\op X_2, U_1\op U_2)$ is a unitary asymptote of $T_1\op T_2$.
\end{cor}

\begin{cor}
\label{restriction}
If $(X,U)$ is a unitary asymptote of the power bounded operator $T\in\L(\HH)$ and $\M\in \textup{Lat}\,T$, then $(X|\M, U|\widetilde\M)$ is a unitary asymptote of the restriction $T|\M$, where $\widetilde\M:= \vee_{n\in\N}U^{-n}X\M$.
\end{cor}

If $T\in\L(\HH)$ is a contraction, then the decreasing sequence $\{\|T^nh\|\}_{n=1}^\i$ is convergent, for every $h\in\HH$.
Suprisingly, these sequences are convergent also in the case, when T is a $\rho$-contraction; see \cite{Eckstein} and \cite{Chevreau}.

If $T\in\L(\HH)$ is an arbitrary power bounded operator, then for every $h\in\HH$ we have
\[\limsup_{n\to\i} \|T^nh\| \le M_T \liminf_{n\to\i} \|T^nh\|.\]
In particular, $\inf\{\|T^nh\|: n\in\N\}=0$ implies $\lim_{n\to\i} \|T^nh\|=0$.
In view of the previous observations we may infer the following statement.

\begin{pro}
If $(X,U)$ is a unitary asymptote of the power bounded operator $T\in\L(\HH)$, then the nullspace of $X$ coincides with the hyperinvariant subspace of those vectors, which are stable for $T$:
\[\textup{ker}\,X= \HH_0(T):= \left\{h\in\HH: \lim_{n\to\i} \|T^nh\|=0\right\}.\]
\end{pro}

We conclude this section with norm-estimates concerning the symbolic calculus and the commutant mapping.

\begin{pro}
\label{norm estimates}
Setting a Banach limit $L$, let $(X_L,U_L)$ be the corresponding unitary asymptote of the power bounded operator $T\in\L(\HH)$.
Then
\begin{itemize}
\item[\textup{(i)}] $\|F\|\le \|\Psi_{X_L}(F)\|\le M_T^2 \|F\| \quad \hbox{ for all } \; F\in\{U_L\}'$;
\item[\textup{(ii)}] $M_T^{-2}\|A_LC\| \le \|\g_{X_L}(C)\| \le \|A_LC\| \quad \hbox{ for all } \; C\in\{T\}'.$
\end{itemize}
\end{pro}

\noindent\textbf{Proof.}
For (i) see the proof of Theorem 3 in \cite{KerToepl}.
Statement (ii) is a consequence of (i) and Proposition \ref{commutant mapping}.

\rightline{$\square$}

\section{Lebesgue decomposition}
\label{decomp}

Let $\P(\T)$ denote the algebra of analytic polynomials, restricted to the unit circle $\T=\{\z\in\C: |\z|=1\}$.
$C(\T)$ is the abelian Banach algebra of all continuous functions defined on $\T$.
The \emph{disk algebra} $\A(\T)$ is the norm-closure of $\P(\T)$ in $C(\T)$.
We recall that the dual of $C(\T)$ can be identified with the Banach space $M(\T)$ of all complex Borel measures on $\T$ (which are automatically regular).
Namely, by the Riesz Representation Theorem, the mapping 
\[\Psi\co M(\T) \to C(\T)^\#, \; \m\mapsto \psi_\m, \quad \hbox{ where }\; \psi_\m(f)= \int_\T f \, d\m,\]
is a Banach space isomorophism.

We say that $\m, \nu\in M(\T)$ are \emph{analytically equivalent measures}, in notation: $\m \ae \nu$, if
\[\int_\T u\, d\m = \int_\T u \, d\nu \quad \hbox{ for all } \, u\in\A(\T).\]
It is clear that this is an equivalence relation.
Furthermore, by a well-known theorem of F. \& M. Riesz,  $\m \ae \nu$ holds if and only if $\nu= \m + h\, dm$ with a function $h\in H_0^1$ (see, e.g., \cite{Hof}).
Here and in the sequel $m$ denotes the normaléized Lebesgue measure on $\T$, and $H_0^1= \{f\in H^1: \int_\T f \, dm= 0\}$.
We note that if the measures $\m, \nu$ are singular (with respect to $m$), then $\m\ae\nu$ holds exactly when $\m=\nu$.

For any $T\in\L(\HH)$, the \emph{polynomial calculus}
\[\FF_{T,0}\co \P(\T)\to \L(\HH),\, p\mapsto p(T)\]
is the uniquely determined (unital) algebra-homomorphism, which transforms the identical function $\h(\z)=\z\,\; (\z\in\T)$ into $T$.
We say that $T$ is a \emph{polynomially bounded operator}, if $\FF_{T,0}$ is norm-continuous.
In that case $\FF_{T,0}$ can be continuously extended to a bounded algebra-homomorphism
\[\FF_{T,1}\co \A(\T)\to\L(\HH), \; u\mapsto u(T).\]
Our aim in this section is to decompose $T$ into the direct sum of an absolutely continuous component and a singular component.
We shall follow Mlak's method of using elementary measures.

More generally, we consider and arbitrary \emph{bounded, linear transformation} $\FF\co \A(\T) \to \L(\HH)$.
For any $x,y\in\HH, \; \f_{x,y}(C)=\la Cx,y\ra$ defines a bounded linear functional on $\L(\HH)$.
By the Hahn--Banach Theorem $\f_{x,y}\circ\FF$ can be extended to a bounded linear functional on $C(\T)$, even in a norm-preserving way.
The set of \emph{elementary measures} of $\FF$ at $x,y$ is defined by
\[M(\FF,x,y):= \left\{\m\in M(\T): \psi_\m|\A(\T)=\f_{x,y}\circ\FF\right\}.\]
Therefore, $\m\in M(\FF,x,y)$ precisely when
\[\la \FF(u)x,y\ra = \int_\T u \, d\m \quad \hbox{ holds for all }\; u\in \A(\T).\]
We say that $\FF$ is \emph{absolutely continuous} (a.c.), if for every $x,y\in\HH$ the measures in $M(\FF,x,y)$ are a.c.\ (with respect to $m$).
On the other hand, $\FF$ is called \emph{singular}, if for every $x,y\in\HH$ the set $M(\FF,x,y)$ contains a measure, which is singular (with respect to $m$).

It is well-known that every measure $\m\in M(\T)$ can be uniquely decomposed into the sum $\m=\m^a + \m^s$, where $\m^a\in M(\T)$ is a.c.\ and $\m^s\in M(\T)$ is singular.
By the next proposition such Lebesgue-type decomposition is valid also for $\FF$.

\begin{pro}
\label{linear decomp}
If $\,\FF\co \A(\T)\to \L(\HH)$ is a bounded linear transformation, then
\begin{itemize}
\item[\textup{(i)}] there exists a unique decomposition $\FF= \FF^a + \FF^s$, where the bounded linear mapping $\,\FF^a\co \A(\T)\to \L(\HH)$ is a.c., and the bounded linear mapping $\FF^s\co \A(\T)\to\L(\HH)$ is singular;
\item[\textup{(ii)}] $\|\FF^a\|\le \|\FF\|$ and $\|\FF^s\|\le \|\FF\|$;
\item[\textup{(iii)}] for every $x,y\in\HH$, we have
$M(\FF^a,x,y)= \{\m^a: \m\in M(\FF,x,y)\}$  and $M(\FF^s,x,y)= \{\m^s+h\, dm: \m\in M(\FF,x,y), h\in H^1_0\}.$
\end{itemize}
\end{pro}

\noindent
\textbf{Proof.}
For every $x,y\in\HH$, the singular component of the measures in $M(\FF,x,y)$ is uniquely determined, denoted by $\m^s_{\FF,x,y}$.
Taking a measure $\m\in M(\FF,x,y)$ with minimal norm, we obtain
\[\|\m^s_{\FF,x,y}\|\le  \|\m\| = \|\f_{x,y}\circ\FF\| \le \|\FF\| \cdot \|x\| \cdot\|y\|.\]
The mapping
\[w^s_\FF\co \HH\times\HH\to  M(\T), \; (x,y)\mapsto \m^s_{\FF,x,y}\]
is sesquilinear (linear in $x$, and conjugate linear in $y$).
For example, let us check the linearity in $x$.
Setting $x_1, x_2, y \in\HH,\; c\in\C$, let us choose $\m_j\in M(\FF,x_j,y_j) \;(j=1,2)$ and $ \m \in M(\FF, c x_1+x_2,y)$.
Then, for every $u\in \A(\T)$, we have
\begin{eqnarray*}
 \int_\T u\, d\m & =  & \la \FF(u)(c x_1+x_2),y\ra = c\la\FF(u)x_1,y\ra + \la\FF(u)x_2,y\ra \\
& = & c \int_\T u \, d\m_1 + \int_\T u \, d\m_2 = \int_\T u \, d(c \m_1+\m_2),
\end{eqnarray*}
whence $\m \ae c\m_1+\m_2$, and so $\m^s= (c\m_1+\m_2)^s=c\m_1^s+\m_2^s$ follows.

For any $g\in C(\T)$, let us consider the bounded linear functional
\[\Lambda_g\co M(\T)\to\C, \; \m\mapsto \int_\T g\, d\m,\]
with $\|\Lambda_g\|=\|g\|$.
The bounded sesquilinear functional
\[w_{\FF,g}^s:= \Lambda_g\circ w_\FF^s\co \HH\times\HH\to \C\]
satisfies the condition $\|w_{\FF,g}^s\|\le \|g\| \|w_\FF^s\|\le \|g\| \|\FF\|$.
There exists a unique operator $\widetilde\FF^s(g)\in\L(\HH)$ such that
\[\left\la\widetilde\FF^s(g)x,y\right\ra = w_{\FF,g}^s(x,y)= \int_\T g \, d\m_{\FF,x,y}^s\]
holds for all $x,y\in\HH$; furthermore $\|\widetilde\FF^s(g)\|= \|w_{\FF,g}^s\| \le \|\FF\|\cdot \|g\|.$
Linearity of the integral in $g$ yields that the mapping $\widetilde\FF^s\co C(\T)\to \L(\HH)$ is linear; we obtain also that $\|\widetilde\FF^s\|\le \|\FF\|$.
Then
\[\FF^s:= \widetilde\FF^s|\A(\T)\co \A(\T)\to \L(\HH)\]
is also a bounded linear transformation, and $\|\FF^s\|\le \|\FF\|$.
Given any $x,y\in\HH$, for every $u\in\A(\T)$, we have
\[\left\la\FF^s(u)x,y\right\ra = \int_\T u\, d\m_{\FF,x,y}^s;\]
hence $M(\FF^s,x,y)= \{\m_{\FF,x,y}^s + h\, dm: h\in H_0^1\}$, and so $\FF^s$ is singular.

Let us consider the bounded linear transformation $\FF^a:= \FF-\FF^s$.
Given any $x,y\in\HH$ and setting $\m\in M(\FF,x,y)$, for every $u\in \A(\T)$, we have
\[\left\la\FF^a(u)x,y\right\ra = \left\la\FF(u)x,y\right\ra - \left\la\FF^s(u)x,y\right\ra = \int_\T u \, d\m - \int_\T u\, d\m^s= \int_\T u\, d(\m-\m^s).\]
Thus $\m^a=\m-\m^s\in M(\FF^a,x,y)$, and so $\FF^a$ is a.c..
Choosing $\m$ so that $\|\m\|= \|\f_{x,y}\circ\FF\|$, we obtain
\[|\left\la \FF^a(u)x,y\right\ra|\le \|u\|\cdot \|\m^a\|\le \|u\|\cdot \|\m\|\le \|u\| \|\FF\| \|x\| \|y\|.\]
Therefore $\|\FF^a(u)\|\le \|\FF\|\cdot \|u\|\; (u\in \A(\T))$, and so $\|\FF^a\|\le \|\FF\|.$

Finally, we turn to the uniqueness of the decomposition.
Suppose that $\FF= \FF_1^a + \FF_1^s = \FF_2^a + \FF_2^s$ are Lebesgue decompositions.
Then $\widetilde\FF= \FF_1^a-\FF_2^a= \FF_2^s - \FF_1^s$ is simultaneously a.c.\ and singular. Thus, the elementary measures of $\widetilde\FF$ are analytically equivalent to 0, and so $\widetilde\FF=0$.

\rightline{$\square$}

\begin{rem}
\label{Banach space setting}
\textup{We note that these considerations can be carried out in the Banach space setting also, with some restrictions.
Let $\X$ be a complex Banach space, and let $\FF\co \A(\T)\to \L(\X)$ be a bounded linear transformation.
For any $x\in\X$ and $y\in\X^\#$, the measure $\m\in M(\T)$ belongs to $M(\FF,x,y)$ if
\[\int_\T u\, d\m= [\FF(u)x,y] \quad \hbox{ for all } \; u\in\A(\T).\]
For any $g\in C(\T)$, we can define the bounded bilinear functional \[w_{\ff,g}^s\co \X\times\X^\#\to\C\]
 as before.
Then a necessary and sufficient condition for the existence of an operator $\widetilde\FF^s(g)\in \L(\X)$ satisfying
\[\left[\widetilde\FF^s(g)x,y\right] = w_{\FF,g}^s(x,y), \quad \hbox{ for all } \; x\in\X, y\in \X^\#,\]
is that $w_{\FF,g}^s(x,y)$ be weak-$*$ continuous in $y$.}
\end{rem}

Now we turn to bounded representations of $\A(\T)$.
We recall that a mapping $\FF\co \A(\T)\to \L(\HH)$ is called a \emph{representation}, if it is a unital algebra-homomorphism.
Norm-continuous representations arise as functional calculae for polynomially bounded operators: $\FF= \FF_{T,1}$ where $T=\FF(\h)$.

Uniqueness of the Lebesgue decomposition of measures implies transformation rules for elementary measures.
Given any $\m, \nu\in M(\T)$, if $\m \ae \nu$ then $\nu= \m + h\, dm$ with some $h\in H_0^1$, and so $\m^a \ae \nu^a$ and $\m^s = \nu^s$.
Furthermore, for any $v\in\A(\T)$, the equation $v\, d\m= v\, d\m^a + v\, d\m^s$ yields that $(v\, d\m)^a = v\, d\m^a$ and  $(v\, d\m)^s = v\, d\m^s$.

\begin{lem}
\label{trf rules}
Let $\FF\co \A(\T)\to \L(\HH)$ be a bounded representation.
Setting $x,y\in \HH$ and $v\in \A(\T)$, let $x_v:= \FF(v)x$ and $y_v^*:= \FF(v)^*y$.
If $\m_{x,y}\in M(\FF,x,y), \m_{x_v,y}\in M(\FF,x_v,y)$ and $\m_{x,y_v^*}\in M(\FF,x,y_v^*)$, then
\[\m_{x_vy}^a \ae v\, d\m_{x,y}^a, \quad \m_{x_v,y}^s = v \, d\m_{x,y}^s, \quad \m_{x,y_v^*}^a \ae v\, d\m_{x,y}^a, \quad \m_{x,y_v^*}^s = v\, d\m_{x,y}^s.\]
\end{lem}

\noindent
\textbf{Proof.}
For every $u\in\A(\T)$, we have
\[\int_\T uv\, d\m_{x,y} = \la\FF(uv)x,y\ra = \la\FF(u)\FF(v)x,y\ra= \la\FF(u)x_v,y\ra = \int_\T u \, d\m_{x_v,y},\]
whence $\m_{x_v,y} \ae v\, d\m_{x,y}$, and so $\m_{x_v,y}^a \ae v\, d\m_{x,y}^a$ and $\m_{x_v,y}^s = v\, d\m_{x,y}^s$ follow.
The other two relations can be derived similarly from the equation $\la\FF(uv)x,y\ra = \la\FF(u)x,y_v^*\ra$.

\rightline{$\square$}

\medskip

We say that the \emph{ polynomially bounded operator} $T\in\L(\HH)$ is \emph{absolutely continuous (singular)} if its functional calculus $\FF_{T,1}$ is an absolutely continuous (singular, respectively) representation.

\begin{thm}
\label{repr decomp}
Let $T\in\L(\HH)$ be a polynomially bounded operator, and let us consider the Lebesgue decomposition $\FF= \FF^a + \FF^s$ of the bounded representation $\FF= \FF_{T,1}\co \A(\T)\to \L(\HH)$.
Then
\begin{itemize}
\item[\textup{(i)}] $\FF^a$ and $\FF^s$ are also bounded representations;
\item[\textup{(ii)}] $P^a= \FF^a(1)$ and $P^s= \FF^s(1)$ are complementary projections belonging to the bicommutant $\{T\}''$;
\item[\textup{(iii)}] $\HH= \HH^a\, \dot+\, \HH^s$, where $\HH^a= P^a\HH$ and $\HH^s=P^s\HH$ are hyperinvariant subspaces of $T$;
\item[\textup{(iv)}] $T^a= T|\HH^a$ is an a.c. and $T^s= T|\HH^s$ is a singular polynomially bounded operator;
the Lebesgue decomposition $T= T^a\, \dot+\, T^s$ of $T$ is unique;
furthermore, $\FF^a(u) = \FF_{T^a,1}(u)\, \dot+\, 0$ and $\FF^s(u)= 0\, \dot+ \,\FF_{T^s,1}(u)$ for all $u\in\A(\T)$.
\end{itemize}
\end{thm}

\noindent\textbf{Proof.}
(i): Given $x,y\in\HH$, set $\m_{x,y}\in M(\FF,x,y)$.
Then $\m_{x,y}^a\in M(\FF^a,x,y)$ and $v\, \m_{x,y}^a \ae \m_{x_v,y}^a$ by Proposition \ref{linear decomp} and Lemma \ref{trf rules}.
Hence, for $u,v\in\A(\T)$, we have
\[\left\la\FF^a(uv)x,y\right\ra = \int_\T uv\, d\m_{x,y}^a = \int_\T u\, d\m_{x_v,y}^a = \left\la\FF^a(u)x_v,y\right\ra = \left\la\FF^a(u)\FF(v)x,y\right\ra.\]
Consequently,
\[\FF^a(uv)= \FF^a(u)\FF(v)\quad \hbox{ for all } \; u,v\in\A(\T).\]
Similarly,
\begin{eqnarray*}
\left\la\FF^a(uv)x,y\right\ra & = & \int_\T uv\, d\m_{x,y}^a = \int_\T v \, d\m_{x,y_u^*}^a \\
& = & \left\la\FF^a(v)x,y_u^*\right\ra = \left\la\FF^a(v)x, \FF(u)^*y\right\ra = \left\la\FF(u)\FF^a(v)x,y\right\ra,
\end{eqnarray*}
and so
\[\FF^a(uv) = \FF(u) \FF^a(v).\]
Then
\begin{eqnarray*}
\FF^a(u)\FF^a(v) + \FF^a(u)\FF^s(v) &= &\FF^a(u)\FF(v)= \FF^a(uv)= \FF(u)\FF^a(v)\\
&=& \FF^a(u)\FF^a(v)+\FF^s(u)\FF^a(v),
\end{eqnarray*}
whence
\[\FF^a(u)\FF^s(v)= \FF^s(u)\FF^a(v) \quad \hbox{ for all } \; u,v\in\A(\T).\]
Fix $v\in\A(\T)$, and set $x,y\in\HH$.
Let $x_1:= \FF^s(v)x$ and $x_2:= \FF^a(v)x$.
Then, for every $u\in\A(\T)$, we have
\begin{eqnarray*}
\int_\T u\, d\m_{x_1,y}^a &=& \left\la\FF^a(u)x_1,y\right\ra = \left\la\FF^a(u)\FF^s(v)x,y\right\ra \\
&=& \left\la\FF^s(u)\FF^a(v)x,y\right\ra = \left\la\FF^s(u)x_2,y\right\ra = \int_\T u \, d\m_{x_2,y}^s.
\end{eqnarray*}
Thus $\m_{x_1,y}^a \ae \m_{x_2,y}^s$, and so $\m_{x_2,y}^s=0$, which implies
 that $\left\la\FF^s(u)\FF^a(v)x,y\right\ra = 0$.
Therefore
\[\FF^a(u)\FF^s(v)= \FF^s(u)\FF^a(v)=0 \quad \hbox{ for all } \; u,v\in\A(\T).\]
Now
\[\FF^a(uv)= \FF^a(u)\FF(v)= \FF^a(u)\FF^a(v)+\FF^a(u)\FF^s(v)= \FF^a(u)\FF^a(v),\]
and so
\begin{eqnarray*}
\FF^s(uv) &=& \FF(uv)-\FF^a(uv)= \left(\FF^a(u)+\FF^s(u)\right)\left(\FF^a(v)+\FF^s(v)\right) - \FF^a(u)\FF^a(v) \\
&=& \FF^s(u)\FF^s(v).
\end{eqnarray*}
Therefore, $\FF^a$ and $\FF^s$ are bounded representations.

\medskip

(ii): It is clear that $(P^a)^2 = (\FF^a(1))^2 = \FF^a(1^2)= P^a$, and similarly $(P^s)^2=P^s$.
Moreover, $P^aP^s= \FF^a(1)\FF^s(1)=0= \FF^s(1)\FF^a(1)= P^sP^a$ and $P^a+P^s= \FF^a(1)+\FF^s(1)= \FF(1)=I$.

Let $C\in\{T\}'$ be arbitrary.
Then $C\FF(u)= \FF(u)C$ for any $u\in\A(\T)$.
Given $x,y\in\HH$, set $\m_{Cx,y}\in M(\FF,Cx,y)$ and $\m_{x,C^*y}\in M(\FF,x,C^*y)$.
Since, for every $u\in\A(\T)$, we have
\[\int_\T u\, d\m_{Cx,y} = \left\la\FF(u)Cx,y\right\ra = \left\la C\FF(u)x,y\right\ra = \left\la\FF(u)x,C^*y\right\ra = \int_\T u \, d\m_{x,C^*y},\]
it follows that $\m_{Cx,y} \ae \m_{x,C^*y}$, and so $\m_{Cx,y}^a \ae \m_{x,C^*y}^a$ and $\m_{Cx,y}^s = \m_{x,C^*y}^s$.
Then, for every $u\in\A(\T)$,
\begin{eqnarray*}
\left\la\FF^a(u)Cx,y\right\ra &=& \int_\T u \, d\m_{Cx,y}^a = \int_\T u \, d\m_{x,C^*y}^a \\
&=& \left\la\FF^a(u)x,C^*y\right\ra = \left\la C\FF^a(u)x,y\right\ra \quad \hbox{ for all } \; x,y\in\HH.
\end{eqnarray*}
Thus $\FF^a(u)C=C\FF^a(u)$, whence $\FF^s(u)C=C\FF^s(u)$ follows.
Therefore
\[\FF^a(u), \FF^s(u)\in\{T\}'' \quad \hbox{ for all } \; u\in\A(\T).\]
In particular, $P^a, P^s \in \{T\}''$.

\medskip

(iii): Since $P^a, P^s\in\{T\}''$, it follows that the subspaces $\HH^a= P^a\HH$ and $\HH^s= P^s\HH$ are hyperinvariant for $T$.
Furthermore, $P^aP^s=P^sP^a=0$ and $P^a+P^s=I$ imply that $\HH^a\, \dot+\, \HH^s =\HH$.

\medskip
(iv): Let us consider the decomposition $T= T^a\, \dot+\, T^s$, where $T^a:= T|\HH^a$ and $T^s:= T|\HH^s$.
For any $u\in\A(\T)$,  we have $\FF^a(u)P^a= P^a\FF^a(u)= \FF^a(1)\FF^a(u)= \FF^a(u)$, whence $\FF^a(u)= (\FF^a(u)|\HH^a)\, \dot+\, 0$ follows.
The bounded representation
\[\FF_0^a\co \A(\T)\to \L(\HH^a), \; \FF_0^a(u):= \FF^a(u)|\HH^a\]
satisfies the condition $\FF_0^a(\h)=T^a=\FF_{T^a,1}(\h)$, and so $\FF_0^a=\FF_{T^a,1}$.
Thus $T^a$ is an a.c.\ polynomially bounded operator.
Similarly, $\FF^s(u)= 0\, \dot+\, (\FF^s(u)|\HH^s)$ for $u\in\A(\T)$, and $\FF_{T^s,1}$ coincides with the singular representation 
\[\FF_0^s\co \A(\T)\to \L(\HH^s), \; \FF_0^s(u):= \FF^s(u)|\HH^s.\]
Thus $T^s$ is a singular polynomially bounded operator.

As for uniqueness, let us suppose that the decomposition $\HH= \HH_1\, \dot+\, \HH_2$ is reducing for $T$, the restriction $T_1=T|\HH_1$ is a.c., and $T_2=T|\HH_2$ is singular.
For every $u\in\A(\T)$, we have $\FF(u)=\FF_{T,1}(u)= \FF_{T_1,1}(u)\, \dot+\, \FF_{T_2,1}(u)$.
Let us consider also the bounded representations $\FF_1, \FF_2\co \A(\T)\to\L(\HH)$ defined by $\FF_1(u):= \FF_{T_1,1}(u)\, \dot+\, 0$ and $\FF_2(u):= 0\, \dot+\, \FF_{T_2,1}(u)$.
Since $T_1$ is a.c., it follows that $\FF_1$ is a.c.; similarly, singularity of $T_2$ implies that $\FF_2$ is singular.
Taking into account that $\FF=\FF_1+\FF_2$, the uniqueness part of Proposition \ref{linear decomp} implies that $\FF_1=\FF^a$ and $\FF_2=\FF^s$.
Thus $\HH_1=\FF_1(1)\HH= \FF^a(1)\HH=\HH^a$ and $\HH_2=\FF_2(1)\HH= \FF^s(1)\HH=\HH^s$.

\rightline{$\square$}

\medskip

Concluding this section we examine the adjoint of a polynomially bounded operator $T\in\L(\HH)$.
For any function $f\co \T\to\C$, let $\widetilde f(\z):= \ol{f(\ol\z)} \; (\z\in\T)$.
For a polynomial $p(\z)= \sum_{n=0}^N c_n\z^n$, we have $\widetilde p(\z)= \sum_{n=0}^N \ol c_n\z^n$, and so $p(T^*)= \widetilde p(T)^*$.
Since $\|p(T^*)\|= \| \widetilde p(T)\|$ and $\|p\|= \|\widetilde p\|$, it follows that $T^*$ is also polynomially bounded, and $\|\FF_{T^*,0}\|= \|\FF_{T,0}\|$.
Taking uniform limits of polynomials, we obtain that for any $u\in\A(\T), \; \widetilde u\in\A(\T)$ and
\[u(T^*)= \widetilde u(T)^*.\]
For any measure $\m\in M(\T), \; \widetilde \m\in M(\T)$ is defined by $\widetilde\m(\o):= \ol{\m(\ol\o)}$, upper bar meaning complex conjugation.
We shall use the notation $M(T,x,y)=M(\FF_{T,1},x,y)$.  

\begin{pro}
\label{adjoint}
If $\, T\in\L(\HH)$ is polynomially bounded, then $T^*$ is also polynomially bounded, and
\begin{itemize}
\item[\textup{(i)}] $M(T^*,x,y)= \{\widetilde \m: \m\in M(T,y,x)\}$ for $ x,y\in\HH$;
\item[\textup{(ii)}] $T$ is a.c.\ if and only if $T^*$ is a.c.;
\item[\textup{(iii)}] $T$ is singular if and only if $T^*$ is singular.
\end{itemize}û
\end{pro}

\noindent\textbf{Proof.}
Given $x,y\in\HH$, set $\m\in M(T,y,x)$ and $\nu\in M(T^*,x,y)$.
For every $u\in \A(\T)$, we have
\begin{eqnarray*}
\int_\T u\, d\nu &=& \la u(T^*)x,y\ra = \la \widetilde u(T)^*x,y\ra= \ol{\la\widetilde u(T)y,x\ra} = \ol{\int_\T \widetilde u \, d\m} \\
&=& \int_\T u(\ol\z)\, d\widetilde\m(\ol\z) = \int_\T u \, d\widetilde\m,
\end{eqnarray*}
whence $\nu \ae \widetilde\m$ follows.
Since $\widetilde h\in H_0^1$, for any $h\in H_0^1$, we obtain (i).
It is easy to see that $\widetilde\m$ is a.c.\ (singular) if and only if $\mu$ is a.c.\ (singular, respectively), which shows the validity of (ii) and (iii).

\rightline{$\square$}

\section{Intertwining relations}
\label{intertw}

Let $T_1\in\L(\HH_1)$ and $T_2\in\L(\HH_2)$ be polynomially bounded operators, and let us assume that $Q\in \I(T_1,T_2)$.
Given $x\in\HH_1$ and $y\in\HH_2$, set $\m\in M(T_2,Qx,y)$ and $\nu\in M(T_1,x,Q^*y)$.
For every $u\in\A(\T)$, we have
\[\int_\T u\, d\m = \la u(T_2)Qx,y\ra = \la u(T_1)x,Q^*y\ra = \int_\T u \, d\nu,\]
hence $\m \ae \nu$, and so
\[M(T_2,Qx,y) = M(T_1,x,Q^*y).\]

\begin{pro}
\label{zero intw}
Let $T_1\in\L(\HH_1)$ be an a.c.\ polynomially bounded operator, and let $T_2\in\L(\HH_2)$ be a singular polynomially bounded operator.
Then $\I(T_1,T_2)= \{0\}$ and $\,\I(T_2,T_1)= \{0\}.$
\end{pro}

\noindent\textbf{Proof.}
Let $Q\in\I(T_1,T_2)$ be arbitrary.
Given $x\in\HH_1$ and $y\in\HH_2$, and setting $\m\in M(T_2,Qx,y)$ and $\nu\in M(T_1,x,Q^*y)$, we obtain from the previous discussion that $\m \ae \nu$.
Since $\m$ is singular and $\nu$ is a.c., it follows that $\m=0$.
Then $\la Qx,y\ra = \int_\T 1\, d\m= 0$; and since $x\in\HH_1, y\in\HH_2$ were arbitrary, we conclude that $Q=0$.
Turning to adjoints and applying Proposition \ref{adjoint}, we obtain that $\I(T_2,T_1)= \{0\}$ also holds.

\rightline{$\square$}

\medskip

We recall that $T_1$ is a \emph{quasiaffine transform} of $T_2$, in notation: $T_1 \prec T_2$, if $\,\I(T_1,T_2)$ contains a quasiaffinity, i.e., an injective transformation with dense range.

\begin{pro}
\label{quasiaffinity}
Let $T_1\in\L(\HH_1)$ and $T_2\in\L(\HH_2)$ be polynomially bounded operators, and let us assume that $T_1 \prec T_2$.
Then
\begin{itemize}
\item[\textup{(i)}] $T_1$ is a.c.\ if and only if $T_2$ is a.c.;
\item[\textup{(ii)}] $T_1$ is singular if and only if $T_2$ is singular.
\end{itemize}
\end{pro}

\noindent\textbf{Proof.}
Let us assume that $T_1$ is a.c., and let us consider the Lebesgue decomposition $T_2= T_2^a\, \dot+\, T_2^s$ of $T_2$.
There exists a quasiaffinity $Q\in\I(T_1,T_2)$.
Let $P_2^s\in \L(\HH_2)$ be the (oblique) projection onto $\HH_2^s$.
Since the transformation $Q^s\in \L(\HH_1,\HH_2^s)$, defined by $Q^sx:= P_2^sQx\; (x\in\HH_1)$, intertwines $T_1$ with $T_2^s$, we infer by Proposition \ref{zero intw} that $Q^s=0$.
Thus $Q\HH_1\subset \HH_2^a$, and since $Q$ has dense range, we obtain that $\HH_2^s= \{0\}$.
Therefore $T_2$ is a.c..
Conversely, assuming that $T_2$ is a.c., the relation $T_2^*\prec T_1^*$ yields that $T_1^*$ is a.c., and then so is $T_1$ too by Proposition \ref{adjoint}.
Statement (ii) can be proved similarly. 

\rightline{$\square$}

\medskip

Let $U\in\L(\K)$ be a unitary operator.
Relying on the Gelfand transform of the abelian $C^*$-algebra generated by $U$, it can be shown that there exists a uniquely determined isometric $*$-representation $\FF\co C(\s(U))\to \L(\K)$ such that $\FF(\h)=U$.
Here $\s(U)\subset\T$ is the spectrum of $U$.
This $\FF$ induces a contractive $*$-representation $\wt\FF\co C(\T)\to\L(\K)$, defined by $f(U)= \wt\FF(f):= \FF(f|\s(U))$.
Since $\wt\FF|\P(\T)= \FF_{U,0}$, it follows that $U$ is polynomially bounded with $\|\FF_{U,0}\|=1$;
furthermore $\wt\FF|\A(\T)= \FF_{U,1}$.
It is known also that $\wt\FF$ can be represented by integration with respect to a uniquely determined spectral measure $E\co \B_\T \to\P(\K)$; see, e.g., \cite{Conway}.
(Here $\B_\T$ denotes the $\s$-algebra of Borel sets on $\T$, and $\P(\K)$ stands for the set of orthogonal projections on $\K$.)
Namely, for every $f\in C(\T)$,
\[\la f(U)x,y\ra = \int_\T f \, dE_{x,y} \quad (x,y\in\K),\]
where $E_{x,y}\in M(\T)$ is the localization of $E$ to $x,y$, defined by $E_{x,y}(\o) := \la E(\o)x,y\ra \; (\o\in\B_\T).$
Therefore
\[M(U,x,y)= \{ E_{x,y}+ h\, dm: h\in H_0^1\}.\]
We recall that the unitary operator $U$ is called a.c.\ (singular, resp.), if $E_{x,y}$ is a.c.\ (singular, resp.) for every $x,y\in\K$.
The previous relation shows that these properties coincide with the corresponding properties considering $U$ as a polynomially bounded operator.

Let us consider the Lebesgue decomposition $U= U^a\, \dot+\, U^s\in \L(\K= \K^a\, \dot+\, \K^s)$.
Since $\K^a\in \hbox{ Hlat}\, U$ and $U^*=U^{-1}\in\{U\}'$, it follows that $\K^a$ is reducing for $U$, and so the orthogonal projection $Q^a\in\L(\K)$ onto $\K^a$ commutes with $U$.
Therefore $Q^a|\K^s\in \I(U^s,U^a)$, whence $Q^a|\K^s=0$ follows by Proposition \ref{zero intw}.
We conclude that $\K^a$ is orthogonal to $\K^s$.
(For a discussion of Lebesgue decomposition of unitaries, based on spectral measures, see \cite{Halmos}.)

Let us consider now a contraction $T\in\L(\HH): \|T\|\le1$.
It is known that $T$ can be decomposed into the orthogonal sum $T= T_c \op U \in \L(\HH= \HH_c \op \HH_u)$ of a completely non-unitary contraction $T_c$ and a unitary operator $U$.
By Sz.-Nagy's Dilation Theorem $T_c$ has a minimal unitary dilation $W$, acting on a Hilbert space $\K$, containing $\HH_c$.
We recall that $W$ is an a.c.\ unitary operator, $T_c^n= P_cW^n|\HH_c$ for every $n\in\Z_+$, and $\K= \vee_{n=1}^\i W^{-n}\HH_c$.
(Here $P_c\in\L(\K)$ is the orthogonal projection onto $\HH_c$.)
In connection with the theory of contractions we refer to \cite{NFBK}.
For every $p\in\P(\T)$, we have
\[\|p(T_c)\|= \|P_cp(W)|\HH_c\|\le \|p(W)\|\le \|p\|,\]
what is called the \emph{von Neumann inequality}, and so $T_c$ is polynomially bounded.
Since, for any $x,y\in\HH_c, \; \la u(T_c)x,y\ra = \la P_cu(W)x,y\ra = \la u(W)x,y\ra$ holds for every $u\in\A(\T)$, we infer that
\[M(T_c,x,y)= M(W,x,y).\]
Therefore $T_c$ is a.c..
Considering the decomposition $T= T_c \op U^a \op U^s$, we can see that the a.c.\ component of $T$ is $T^a= T_c \op U^a$, while the singular part is $T^s=U^s$.

If $T\in\L(\HH)$ is a polynomially bounded operator, $U\in\L(\K)$ is a singular unitary operator, and $T \prec U$, then $T$ is singular by Proposition \ref{quasiaffinity}.
The following theorem, which in different forms can be found in \cite{Ando} and \cite{Mlak algebraic}, states that there must be a much stronger relation between $T$ and $U$ in that case.

\begin{thm}
\label{singular}
If $\, T\in\L(\HH)$ is a singular polynomially bounded operator, then
\begin{itemize}
\item[\textup{(i)}] the funcional calculus $\FF_{T,1}\co \A(\T) \to \L(\HH)$ can be extended to a bounded representation $\wt\FF_{T,1}\co C(\T)\to \L(\HH)$;
\item[\textup{(ii)}] the operator $X\in\I(T,U)$ is invertible, whenever $(X,U)$ is a unitary asymptote of $T$; and so $T$ is similar to the singular unitary operator $U$.
\end{itemize}
\end{thm}

\noindent\textbf{Proof.}
(i): We repeat the procedure carried out in the first part of the proof of Proposition \ref{linear decomp}.
We shall write $\FF=\FF_{T,1}$ and $\wt\FF=\wt\FF_{T,1}$ for short.
For any $x,y\in\HH$, there exists a unique singular measure $\m_{x,y}$ in $M(T,x,y)$.
Let us consider the bounded sesquilinear mapping
\[w_T\co \HH\times\HH \to M(\T), \; w_T(x,y)= \m_{x,y}.\]
For any $g\in C(\T), \; \Lambda_g\co M(\T)\to\C$ is defined by $\Lambda_g(\m)= \int_\T g \, d\m$.
The bounded sesquilinear functional
\[w_{T,g}:= \Lambda_g\circ w_T\co \HH\times\HH\to \C\]
uniquely determines an operator $\wt\FF(g)\in \L(\HH)$:
\[\la\wt\FF(g)x,y\ra = w_{T,g}(x,y)= \int_\T g \, d\m_{x,y}\quad (x,y\in\HH).\]
We know that $\|\wt\FF(g)\|\le \|\FF\| \|g\|$.
It is clear also that $\wt\FF$ is linear.
We have to show yet that $\wt\FF$ is multiplicative.

For $v\in\A(\T)$, we have $v\, d\m_{x,y}= \m_{x_v,y}$ by Lemma \ref{trf rules}.
Thus, for any $f\in C(\T)$,
\[\la\wt\FF(fv)x,y\ra = \int_\T fv\, d\m_{x,y}= \int_\T f\, d\m_{x_v,y} = \la\wt\FF(f)x_v,y\ra = \la\wt\FF(f)\FF(v)x,y\ra,\]
whence $\wt\FF(fv)= \wt\FF(f)\FF(v)$ follows.
Setting $y_f^*= \wt\FF(f)^*y$, for every $v\in\A(\T)$,
\[\int_\T vf\, d\m_{x,y} = \la\wt\FF(fv)x,y\ra = \la\wt\FF(f)\FF(v)x,y\ra = \la \FF(v)x,y_f^*\ra = \int_\T v \, d\m_{x,y_f^*}.\]
Thus, the singular measures $f\, d\m_{x,y}$ and $\m_{x,y_f^*}$ are analytically equivalent, and so they must coincide: $f\, d\m_{x,y}= \m_{x,y_f^*}$.
Taking any $g\in C(\T)$,
\[\la\wt\FF(fg)x,y\ra = \int_\T fg \, d\m_{x,y}= \int_\T g\, d\m_{x,y_f^*} = \la\wt\FF(g)x, y_f^*\ra = \la\wt\FF(f)\wt\FF(g)x,y\ra\]
holds for every $x,y\in\HH$. 
Consequently, $\wt\FF(fg)= \wt\FF(f)\wt\FF(g)$.

\medskip

(ii):
By the multiplicativity of $\wt\FF$ we infer that $T$ is invertible and $\wt\FF(\h^n)=T^n$ holds, for every $n\in\Z$.
Since $\|T^n\|\le \|\wt\FF\| = \|\FF\|$ for all $n\in\Z$, a well-known theorem of Sz.-Nagy yields that $T$ is similar to a unitary operator $V$; see \cite{Nagy47}.
Let $Q\in \I(T,V)$ be invertible.
The polynomially bounded operator $T$ is necessarily power bounded, and so it has a unitary asymptote $(X,U)$.
There exists a unique $Y\in\I(U,V)$ such that $Q=YX$.
Since $Q$ is invertible, it follows that $X$ is bounded from below.
The isometry $U|X$
is singular by Proposition \ref{quasiaffinity}. 
Considering its Wold decomposition, we can see that $U|X$ is unitary.
Hence $X\HH$ is reducing for $U$, and so the minimality of $(X,U)$ yields that $X$ is a surjection.
Therefore $X$ is invertible, and then $Y$ must be invertible too.
Consequently, the operators $T, U$ and $V$ are similar to each other.

\rightline{$\square$}

\medskip

We conclude this section with a discussion of further intertwining relations.
Since Hilbert space contractions have a rich theory (see \cite{NFBK}), it would be important to know how more general operators can be related to contractions.
Answering a question posed by Sz.-Nagy, Foguel gave examples for power bounded operators which are not similar to contractions; see \cite{Foguel}.
Pisier answered negatively also a more delicate question of Halmos, showing that not every polynomially bounded operator is similar to a contraction; see \cite{Pisier}.
M\"uller and Tomilov proved, giving negative answer for a question of the author posed in \cite{KerInd89}, that there are power bounded operators of class $C_{11}$ which are not similar to contractions; see \cite{Muller}.
We recall that $T\in C_{11}$ means $\HH_0(T)=\HH_0(T^*)=\{0\}$, and in that case $T$ is \emph{quasisimilar} to its unitary asymptote $U: T \sim U$, that is $T \prec U$ and $U \prec T$ hold simultaneously.

\begin{que} If the polynomially bounded operator $T$ is of class $C_{11}$, is $T$ similar to a contraction?
\end{que}

For additional conditions under which a power bounded operator is similar to a contraction see \cite{Gamal}.
We note yet that by a fundamental characterization due to Paulsen, exactly those operators are similar to contractions which are completely bounded; see \cite{Paulsen}.

Quasisimilarity is a much weaker relation than similarity, but it preserves also numerous properties and plays important role in the classification of operators; see, e.g., \cite{Bercovici} and \cite {Davidson}.
M\"uller and Tomilov showed that a power bounded operator $T$ is not necessarily quasisimilar to a contraction \cite{Muller}.
On the other hand, Bercovici and Prunaru proved that if $T$ is a polynomially bounded operator, then there exist contractions $T_1$ and $T_2$ such that $T_1 \prec T \prec T_2$; see \cite{BP}.
As far as we know, the following questions are still open.

\begin{que}
Is every polynomially bounded operator quasisimilar to a contraction?
\end{que}

\begin{que}
If the power bounded operator $T$ is quasisimilar to a singular unitary operator $V$, does it follow that $T$ is similar to $V$?
\end{que}

The latter question is connected with Theorem \ref{singular}, and was posed in \cite{KerInd89}.
Partial answer for it can be found in \cite{Gamal sing}.

\section{$H^\i$-functional calculus}
\label{calculus}

The Hardy class $H^\i$ is the weak-$*$ closed subalgebra
of $L^\i(\T)$, consisting of those functions $f$ whose Fourier coefficients satisfy the condition $\widehat f(-n)=0$ for $n\in\N$.
$H^\i$ can be identified as the dual of the Banach space $L^1(\T)/H_0^1$.
For a detailed study of this class see \cite{Hof}.
Sz.-Nagy and Foias introduced $H^\i$-functional calculus for a.c.\ contractions, and thoroughly exploited its properties in their theory of contractions; see \cite{NFBK} and \cite{Bercovici}.
In Chapter 2 of \cite{Challendar} $H^\i$-functional calculus is defined for a polynomially bounded operator $T$ acting on a complex Banach space $\X$, which is stable, that is $\lim_{n\to\i} \|T^nx\|=0$ holds for avery $x\in\X$.
In \cite{Mlak Acta} and \cite{Mlak Oberwolfach} Mlak considered representations of general function algebras.
Here we follow Mlak's method of elementary measures, concentrating on $H^\i$, and providing detailed study of the calculus in this case.
We note that in \cite{KerToepl} our approach was based on the unitary asymptote in the $C_{1\cdot}$-case, that is when $\HH_0(T)=\{0\}$.

We recall that $\L(\HH)$ is the dual of the Banach space ${\cal C}_1(\HH)$ of trace class operators.
Namely, the mapping $\Lambda\co \L(\HH)\to {\cal C}_1(\HH)^\#$, defined by $[A,\Lambda(C)]= \hbox{tr}(AC)$, is a Banach space isomorphism; see \cite{Schatten}.

We say that the operator $T\in\L(\HH)$ \emph{admits an $H^\i$-functional calculus}, if there exists a weak-$*$ continuous representation $\FF_{T,2}\co H^\i \to \L(\HH)$ such that $\FF_{T,2}(\h)=T$.
In that case we use the notation $\FF_{T,2}(f) = f(T)\; (f\in H^\i)$.
In the following proposition we collect some basic facts about this calculus.

\begin{pro}
\label{basic facts}
If $T\in\L(\HH)$ admits an $H^\i$-functional calculus, then
\begin{itemize}
\item[\textup{(i)}] $\FF_{T,2}$ is uniquely determined;
\item[\textup{(ii)}] $\FF_{T,2}$ is norm-continuous;
\item[\textup{(iii)}] $T$ is polynomially bounded, $\FF_{T,2}|\A(\T)= \FF_{T,1}$, and $\|\FF_{T,2}\| = \|\FF_{T,1}\|$;
\item[\textup{(iv)}]  $H^\i(T):= \textup{ ran}\,\FF_{T,2}\subset \{T\}''$;
\item[\textup{(v)}]  $\textup{ Lat}\,T = \textup{ Lat}\,H^\i(T)$;
\item[\textup{(vi)}] for every $\M\in \textup{ Lat}\,T, \; T|\M$ admits an $H^\i$-functional calculus, and $f(T|\M)= f(T)|\M\; (f\in H^\i).$
\end{itemize}
\end{pro}

\noindent\textbf{Proof.}
(i): Let us assume that $\FF_j\co H^\i\to \L(\HH)$ is a weak-$*$ continuous representation, for $j=1,2$, and that $\FF_1(\h)=\FF_2(\h)$ holds.
Then $\FF_1(p)= \FF_2(p)$ is true, for every polynomial $p\in \P(\T)$.
Taking an arbitrary $f\in H^\i$, the Cesaro means 
\[\s_{f,n} = \sum_{k=0}^n \left(1-\frac{k}{n+1}\right)\widehat f(k)\h^k\]
of the Fourier series of $f$ converge to $f$ in the weak-$*$ topology, as $n\to\i$; see, e.g., Chapter 2 in \cite{Hof}.
By continuity we obtain that $\FF_1(f)= \FF_2(f)$.

(ii): For short we write $\FF= \FF_{T,2}$.
The unit ball $(H^\i)_1:= \{f\in H^\i: \|f\|\le 1\}$ is weak-$*$ compact by the Banach--Alaoglu Theorem.
Since $\FF$ is weak-$*$ continuous, it follows that $\FF(H^\i)_1$ is weak-$*$ compact in $\L(\HH)$.
Thus $\{[A,\FF(f)] : f\in (H^\i)_1\}$ is compact and so bounded in $\C$, for every $A\in {\cal C}_1(\HH)$.
Then the Uniform Boundedness Principle yields that $\FF(H^\i)_1$ is bounded in $\L(\HH)$, which means the norm-continuity of $\FF$.

(iii): Since $\FF|\P(\T)= \FF_{T,0}$, we infer from (ii) that $\FF_{T,0}$ is bounded: $\|\FF_{T,0}\|\le \|\FF\|$.
Hence $T$ is polynomially bounded, and $\FF|\A(\T)= \FF_{T,1}$.
Furthermore, for any $f\in H^\i$, we have $\|\FF(\s_{f,n})\|\le \|\FF_{T,0}\|\cdot \|f\|\; (n\in\N).$
Taking into account that $w^*$-$\lim_n\FF(\s_{f,n})= \FF(f)$, we conlude that $\|\FF(f)\|\le \|\FF_{T,0}\|\cdot \|f\|$.
Therefore $\|\FF\|= \|\FF_{T,0}\|= \|\FF_{T,1}\|.$

The statements (iv), (v), (vi) can be easily derived from the fact that $H^\i$ is the (sequentially) weak-$*$ closure of $\P(\T)$.

\rightline{$\square$}

\medskip
The following lemma shows that it is enough to check seemingly weaker conditions in order to prove that a mapping is an $H^\i$-functional calculus.

\begin{lem}
\label{seemingly weaker}
If $\,\FF\co H^\i\to \L(\HH)$ is a linear mapping, $\FF|\A(\T)$ is a representation, and $\lim_n \la\FF(f_n)x,y\ra=0$ for every $x,y\in\HH$ whenever $w^*$-$\lim_n f_n=0$, then $\FF$ is an $H^\i$-functional calculus.
\end{lem}

\noindent\textbf{Proof.}
Let us assume that the sequence $\{f_n\}_{n=1}^\i$ in $H^\i$ converges to zero in the weak-$*$ topology.
By the assumption, the operators $\{\FF(f_n)\}_{n=1}^\i$ converge to zero in the weak operator topology (wot).
The Uniform Boundedness Principle yields that $M= \sup\{\|\FF(f_n)\|: n\in\N\}<\i$.
Taking an arbitrary operator $A\in {\cal C}_1(\HH)$, let us consider the Hilbert--Schmidt decomposition $A= \sum_k s_k\, x_k \otimes y_k$.
Here $\{x_k\}_k$ and $\{y_k\}_k$ form orthonormal systems, $s_k\ge0$ for all $k$, and $\sum_k s_k= \|A\|_1<\i$.
These conditions imply that
\[\lim_{n\to\i} \textup{tr}(\FF(f_n)A) = \lim_{n\to\i} \sum_k s_k \la\FF(f_n)x_k,y_k\ra=0.\]
Thus the functional $\Lambda_A\circ\FF$ is sequentially weak-$*$ continuous, where the weak-$*$ continuous functional $\Lambda_A\co \L(\HH)\to\C$ is defined by $\Lambda_A(C)= \textup{ tr}(CA)$.
By the Krein--Smulian Theorem $\Lambda_A\circ\FF$ is weak-$*$ continuous; see Corollary V.12.8 in \cite{Conway}.
Therefore, $\FF$ is also weak-$*$ continuous; see Proposition 1.3.2 in \cite{Kadison}.

Setting any $f,g\in H^\i$, we know that $\FF(\s_{f,n}\s_{g,k})= \FF(\s_{f,n})\FF(\s_{g,k})$ holds for every $n,k\in \N$.
Since $w^*$-$\lim_n\s_{f,n}=f$, we infer by the sequentially weak-$*$--wot continuity of $\FF$ that $\FF(f\s_{g,k})= \FF(f)\FF(\s_{g,k})$.
Tending now $k$ to infinity, we conclude that $\FF(fg)= \FF(f)\FF(g)$.

\rightline{$\square$}

\begin{thm}
\label{existence of functional calculus}
The operator $T\in\L(\HH)$ admits an $H^\i$-functional calculus, if and only if $\, T$ is an absolutely continuous polynomially bounded operator.
\end{thm}

\noindent\textbf{Proof.}
Let us assume that $T$ admits an $H^\i$-functional calculus.
Then $T$ is necessarily polynomially bounded by Proposition \ref{basic facts}.
Suppose that $T$ is not a.c., and consider the Lebesgue decomposition $T= T^a  \,\dot+\, T^s\in \L(\HH= \HH^a\, \dot+\, \HH^s)$, where $\HH^s\neq \{0\}$. 
In view of Proposition \ref{basic facts}, $T^s$ admits an $H^\i$-functional calculus.
We know by Theorem \ref{singular} that $T^s$ is similar to a singular unitary operator $U\in\L(\K)$, let $Q\in\I(T,U)$ be invertible.
It is clear that $\FF_{U,2}(f):= Q\FF_{T,2}(f)Q^{-1}\; (f\in H^\i)$ defines an $H^\i$-calculus for $U$.
Let $E$ denote the spectral measure of $U$.
Taking a non-zero vector $x\in\K$, let us consider the positive measure $E_{x,x}$. 
We know that $E_{x,x}$ is singular and $E_{x,x}(\T)= \|x\|^2>0$.
In view of regularity we can find a compact set $K\subset\T$ such that $m(K)=0$ and $E_{x,x}(K)>0$.
By a result of Rudin there exists a function $u\in\A(\T)$ such that $u(\z)=1$ for all $\z\in K$, and $|u(\z)|<1$ for all $\z\in \T\setminus K$; see page 81 in \cite{Hof}.
Since $w^*$-$\lim_n u^n=0$, regarding the continuity properties of $\FF_{U,2}$ it follows that $\lim_n \la u^n(U)x,x\ra=0$.
However, 
\[\lim_{n\to\i}\la u^n(U)x,x\ra = \lim_{n\to\i} \int_\T u^n\, dE_{x,x} = E_{x,x}(K)>0.\]
We arrived at a contradiction, and so the polynomially bounded operator $T$ must be absolutely continuous.

Let us assume now that $T$ is an a.c.\ polynomially bounded operator.
Let $f\in H^\i$.
Given any $x,y\in \HH$, let us choose an elementary measure $\m_{x,y}\in M(T,x,y)$.
Being abolutely continuous, $\m_{x,y}= g_{x,y}\, dm$ where $g_{x,y}\in L^1(\T)$; hence the integral $\int_\T f \, d\m_{x,y}= \int_\T fg_{x,y}\, dm$  can be formed.
Any other $\nu\in M(T,x,y)$ is of the form $\nu= \m_{x,y} + h\, dm$ with $h\in H_0^1$, and so $\int_\T f \, d\nu = \int_\T f\, d\m_{x,y} + \int_\T fh\, dm = \int_\T f\, d\m_{x,y}$, because $fh\in H_0^1$.
Let us consider the well-defined mapping
\[w_f\co \HH\times\HH\to \C, \; w_f(x,y)= \int_\T f \, d\m_{x,y}.\]
It is easy to check that $w_f$ is sesquilinear, and $\|w_f\|\le \|f\|\cdot \|\FF_{T,1}\|.$
There exists a unique operator $\FF(f)\in\L(\HH)$ such that
\[\la \FF(f)x,y\ra= w_f(x,y) = \int_\T f\, d\m_{x,y} \quad (x,y\in\HH).\]
It is clear that $\FF\co H^\i\to\L(\HH)$ is linear, and $\FF|\A(\T)= \FF_{T,1}$ is a representation.
If $\{f_n\}_{n=1}^\i$ is a sequence in $H^\i$ converging to 0 in the weak-$*$ topology, then
\[\lim_{n\to\i} \la\FF(f_n)x,y\ra = \lim_{n\to\i} \int_\T f_ng_{x,y}\, dm=0 \quad (x,y\in\HH). \]
Applying Lemma \ref{seemingly weaker} we obtain that $\FF$ is an $H^\i$-functional calculus for $T$.

\rightline{$\square$}

\begin{pro}
\label{function of adjoint}
If $\, T\in\L(\HH)$ admits an $H^\i$-functional calculus, then so does its adjoint $T^*$.
Furthermore, $f(T^*)= \wt f(T)^*$ holds for every $f\in H^\i$.
\end{pro}

\noindent\textbf{Proof.}
We have seen that $T$ is an a.c.\ polynomially bounded operator if and only if so is $T^*$; see Proposition \ref{adjoint}.
Furthermore, for any polynomial $p\in\P(\T)$ we have $p(T^*)= \wt p(T)^*$.
Let $f\in H^\i$ be arbitrary.
Since $w^*$-$\lim_n \s_{f,n} = f$ and $w^*$-$\lim_n \wt\s_{f,n}=\wt f$, it follows that $\{\s_{f,n}(T^*)\}_n$ converges to $f(T^*)$, and $\{\wt\s_{f,n}(T)\}_n$ converges to $\wt f(T)$ in the weak operator topology.
The latter condition implies that $\{\wt\s_{f,n}(T)^*\}_n$ converges to $\wt f(T)^*$ in wot.
Therefore, the equalities $\s_{f,n}(T^*)= \wt\s_{f,n}(T)^*\; (n\in\N)$ yield that $f(T)^*= \wt f(T)^*.$

\rightline{$\square$}

\section{Quasianalytic operators}
\label{quasianalytic}

Let $T\in\L(\HH)$ be an a.c.\ polynomially bounded operator.
Let $(X,U)$ be a unitary asymptote of $T$, where $U\in\L(\K)$.
Taking the Lebesgue decomposition $U= U^a \op U^s\in \L(\K= \K^a\op\K^s)$, and considering the transformation $Y\in\I(T,U^s)$, defined by $Yh=P^sXh$, we infer by Proposition \ref{zero intw} that $Y=0$.
Hence $X\HH\subset\K^a$, and it follows by the minimality of $(X,U)$ that $\K^s=\{0\}$.
Therefore, \emph{$U$ is an a.c.\ unitary operator}.

Let $E\co \B_\T\to\P(\K)$ denote the spectral measure of $U$.
For any $x,y\in \HH$, the localization of $E$ at $Xx, Xy$ is of the form $E_{Xx,Xy}=w_{x,y} \ dm$, where $w_{x,y}\in L^1(\T)$ is the asymptotic density function of $T$ at $x,y$.
We note that
\[E_{Xx,Xy}\in M(T,x,Ay)\]
with $A=X^*X$ (see the beginning of Section \ref{intertw}), and $E_{Xx,Xx}$ is the unique positive measure in $M(T,x,Ax)$.

The measurable set
\[\o(T,x):= \{\z\in\T: w_{x,x}(\z)>0\}\]
is called the  \emph{local residual set} of $T$ at $x$.
Considering the functional model of the unitary operator $U$, we may easily check that the hyperinvariant subspace of $U$ generated by $Xx$ is just the spectral subspace corresponding to $\o(T,x)$:
\[E(\o(T,x))\K = \vee\{FXx: F\in\{U\}'\} \in \textup{ Hlat}\,U.\]
More precisely, $\o(T,x)$ is the smallest measurable set on $\T$ such that $E(\o(T,x))\K$ contains $Xx$, that is $Xx\in E(\a)\K$ implies $m(\o(T,x)\setminus\a)=0$; see, e.g.,   \cite{Conway}.

\begin{pro}
\label{unique}
The local residual set $\o(T,x)$ is independent of the particular choice of the unitary asymptote $(X,U)$.
\end{pro}

\noindent\textbf{Proof.}
Let $(X',U')$ be another unitary asymptote of $T$; $E'$ is the spectral measure of $U'\in\L(\K')$, and $\o'(T,x)$ is the local residual set of $T$ at $x$, defined via $E'$.
We know that there is an invertible transformation $Z\in\I(U,U')$ such that $X'=ZX$.
Let us consider the polar decomposition $Z= W |Z|$, where $|Z|\ge0$ is invertible and $W$ is unitary.
It is easy to verify that $|Z|\in \{U\}'$ and $WU=U'W$.
Since the vector $X'x=ZXx$ belongs to the subspace
\[ZE(\o(T,x))\K= W |Z| E(\o(T,x))\K = W E(\o(T,x)) |Z|\K = E'(\o(T,x)) \K',\]
it follows that $m(\o'(T,x)\setminus \o(T,x))=0$.
Changing the roles of $(X,U)$ and $(X',U')$ we obtain that the symmetric difference $\o(T,x)\,\triangle\,\o'(T,x)$ is of zero Lebesgue measure, and so they can be considered identical.
Note that the Radon--Nikodym derivative $w_{x,x}$, and so $\o(T,x)$ also are determined up to sets of zero Lebesgue measure.

\rightline{$\square$}

\medskip

Let us consider the functional calculus $\FF_{T,2}\co H^\i\to\L(\HH), \; f\mapsto f(T)$ for $T$.
We write $K_T= \|\FF_{T,2}\|$ for short.
We recall that for $f,g\in H^\i, \; f \overset{a}\prec g$ means that $|f(z)|\le |g(z)|$ for all $z\in\D:= \{\z\in\C: |\z|<1\}$.
Then $f=gh$, where $h\in H^\i$ and $\|h\|\le1$.
Thus $\|f(T)x\|= \|h(T)g(T)x\|\le K_T \|g(T)x\| \; (x\in\HH)$, that is $f(T) \overset{a}\prec K_T g(T).$
(If $T$ is a contraction, then $K_T=1$ and so $\FF_{T,2}$ is monotone.)
Let $F=\{f_n\}_{n=1}^\i$ be a decreasing sequence in $H^\i$, i.e.\ $f_{n+1} \overset{a}\prec f_n$ for every $n$.
The measurable limit function $\f_F(\z)= \lim_n |f_n(\z)|$ is defined for almost every $\z\in\T$.
Set
\[N_F := \{\z\in\T: \f_F(\z)>0\}.\]
We say that $F$ is asymptotically non-vanishing on the measurable set $\a\subset\T$, if $m(\a\cap N_F)>0$.

\begin{pro}
\label{converging to zero}
If $\, \inf\{\|f_n(T)x\|: n\in\N\}=0$, then $\lim_{n\to\i}\|f_n(T)x\|=0$.
Furthermore,
\[\HH_0(T,F):= \{x\in\HH: \lim_{n\to\i} \|f_n(T)x\|=0\}\]
is a hyperinvariant subspace of $T$.
\end{pro}

\noindent\textbf{Proof.}
The first statement follows from the condition $f_{n+1}(T) \overset{a}\prec K_T f_n(T)\; (n\in\N)$, while the second one is a consequence of $\sup\{\|f_n(T)\|: n\in\N\}<\i.$

\rightline{$\square$}

\medskip

The following theorem shows that the local residual set is responsible for local stability.

\begin{thm}
\label{local stability}
\begin{itemize}
\item[\textup{(a)}] If $x\in \HH_0(T,F)$, then $m(N_f\cap\o(T,x))=0$.
\item[\textup{(b)}] If $m(N_F\cap \o(T,x))=0$, then there exists a strictly increasing mapping $\tau\co \N\to\N$ such that $x\in\HH_0(T,G)$, where $G=\{\h^{\tau(n)}f_n\}_{n=1}^\i$.
\end{itemize}
\end{thm}

\noindent\textbf{Proof.}
For any $x\in\HH$, we have
\[\|Xf_n(T)x\|^2= \|f_n(U)Xx\|^2 = \int_\T|f_n|^2w_{x,x}\, dm \quad (n\in\N),\]
and the latter integral converges to $\int_\T \f_F^2 w_{x,x} \, dm$ by Lebesgue's Dominated Convergence Theorem, as $n\to\i$.
Hence (a) follows, and we obtain (b) also by applying Theorem \ref{characterization}.
Note that $\f_G= \f_F$.

\rightline{$\square$}

\medskip
The previous theorem tells us that if the decreasing sequence $F=\{f_n\}_{n=1}^\i$ in $H^\i$ is asymptotically non-vanishing on $\o(T,x)$, then $x\not\in \HH_0(T,F)$, that is the vector-sequence $\{f_n(T)x\}_{n=1}^\i$ is asymptotically non-vanishing.
On the other hand, if $F$ is asymptotically vanishing on $\o(T,x)$, then  a modified (strengthened) sequence $\{T^{\tau(n)}f_n(T)x\}_{n=1}^\i$ is asymptotically vanishing (stable) with a suitable mapping $\tau$.

It is known that the measurable sets on $\T$ form a complete lattice, if we disregard of sets of measure zero.
Hence, for any measurable sets $\a, \b$ on $\T$, we shall write $\a\subset\b, \; \a=\b$ or $\a\neq\b$, when $m(\a\setminus\b)=0, \; m(\a\,\triangle\,\b)=0$ or $m(\a\,\triangle\,\b)>0$, respectively.
The following statement can be proved using the Gelfand transform of the abelian Banach algebra $L^\i(\T)$; see Section 11.13 in \cite{Rudin}.
Here we sketch an elementary proof.

\begin{lem}
\label{measure lattice}
For any system of measurable sets $\{\o_i\}_{i\in I}$ on $\T$, there exist a smallest measurable set $\ol\o= \vee_{i\in I}\o_i$ containing all $\o_i$, and a largest measurable set $\underline\o= \wedge_{i\in I}\o_i$ contained in all $\o_i$.
Furthermore, $\ol\o$ and $\underline\o$ are uniquely determined.
\end{lem}

\noindent\textbf{Proof.}
Let $\ol\O$ be the system of those measurable sets, which contain (in the extended sense) every $\o_i\; (i\in I)$.
It is clear that $\ol\O$ is closed under countable intersection.
Set $a= \inf\{m(\o): \o\in\ol\O\}$, and select a sequence $\{\widehat\o_n\}_{n=1}^\i$ in $\ol\O$ so that $\lim_n m(\widehat\o_n)=a$.
Then it is easy to verify that the measurable set $\ol\o= \cap_{n=1}^\i \widehat\o_n$ has the required properties.
The statement about $\underline\o$ can be proved similarly, considering the system $\underline\O$ of those measurable sets, which are contained in every $\o_i\; (i\in I)$, and then taking $b= \sup\{m(\o): \o\in\underline\O\}$.

\rightline{$\square$}

\medskip

The set 
\[\o(T):= \vee\{\o(T,x): x\in\HH\}\]
is called the \emph{residual set} of $T$.
It is clear that $\o(T)$ is the measurable support of the spectral measure $E$ of $U$, that is $E(\a)=0$ exactly when $m(\a\cap\o(T))=0$.
Furthermore, $\o(T)=\emptyset$ means that $\K=\{0\}$, that is $\HH_0(T)=\HH$, in which case $T$ is called of class $C_{0\cdot}$.

The set
\[\pi(T):= \wedge\{\o(T,x): 0\neq x\in\HH\}\]
is called the \emph{quasianalytic spectral set} of $T$.
In view of Theorem \ref{local stability}, this is the largest measurable set with the property that $\HH_0(T,F)=\{0\}$, whenever $F=\{f_n\}_{n=1}^\i$ is a decreasing sequence in $H^\i$, asymptotically non-vanishing on $\pi(T)$.
We can characterize $\pi(T)$ also with more general sequences, as in the next proposition.
(For related results in connection with test sequences of stability, see \cite{KSztest}.)

\begin{pro}
\label{bounded sequences}
If $F=\{f_n\}_{n=1}^\i$ is a bounded sequence in $H^\i$ such that $\limsup_n \|\h_{\pi(T)}f_n\|_2>0$, then the hyperinvariant subspace \[\HH_0(T,F):= \{x\in\HH: \lim_{n\to\i}\|f_n(T)x\|=0\}\]
reduces to zero.
Furthermore, $\pi(T)$ is the largest measurable set having this property.
\end{pro}

\noindent\textbf{Proof.}
Let $F=\{f_n\}_{n=1}^\i$ be a bounded sequence in $H^\i$ such that $c= \limsup_n\int_{\pi(T)}|f_n|^2\, dm>0$.
Set $M= \sup\{\|f_n\|_\i : n\in\N\}< \i$.
Furthermore, given $0\neq x\in\HH$, let $\a(x,N):= \{\z\in\pi(T): w_{x,x}(\z)> N^{-1}\}$ for $N\in\N$.
We can select $N$ so that $M^2 m(\pi(T)\setminus\a(x,N))< c/2$.
For every $n\in\N$, we have
\begin{eqnarray*}
\int_\T |f_n|^2 w_{x,x}\, dm &\ge& N^{-1} \int_{\a(x,N)}|f_n|^2\, dm \\
&=& N^{-1} \left( \int_{\pi(T)}|f_n|^2\, dm - \int_{\pi(T)\setminus\a(x,N)}|f_n|^2\, dm\right)\\
&\ge& N^{-1} \big(\|\h_{\pi(T)}f_n\|^2_2 - M^2 m(\pi(T)\setminus\a(x,N))\big)
\end{eqnarray*}
Hence
\[\limsup_{n\to\i} \|Xf_n(T)x\|^2= \limsup_{n\to\i} \int_\T |f_n|^2 w_{x,x}\, dm \ge N^{-1} c/2>0,\]
and so $x$ is not stable for $F$.
Therefore $\HH_0(T,F)=\{0\}$.

Let $\b$ be a measurable set on $\T$, larger than $\pi(T): \; \b\supset\pi(T)$ and $\b\neq\pi(T)$.
By the definition of $\pi(T)$ we can find a non-zero $x_0\in\HH$ so that $\b$ is not contained in $\o(T,x_0)$, that is $\Delta= \b\setminus\o(T,x_0)$ is of positive measure.
Let $f_0\in H^\i$ be an outer function satisfying the condition $|f_0|= \h_{\Delta} + (1/2) \h_{\T\setminus\Delta}$, and let us consider the decreasing sequence $F_0=\{f_0^n\}_{n=1}^\i$ in $H^\i$, with $\f_{F_0}=\h_{\Delta}$.
In view of Theorem \ref{local stability}, there exists a strictly increasing $\tau\co \N\to\N$ such that $x_0\in\HH_0(T,F)$ with $F= \{\h^{\tau(n)}f_0^n\}_{n=1}^\i$.
It is clear also that
\[\lim_{n\to\i} \|\h_\b \h^{\tau(n)}f_0^n\|_2^2 = m(\Delta)>0,\]
which completes the proof.

\rightline{$\square$}

\medskip
The absolutely continuous polynomially bounded operator $T\in\L(\HH)$ is \emph{quasianalytic}, if $\pi(T)=\o(T)\neq\emptyset$.
We note that $\pi(T)\subset\o(T)$ is evident from the definition.
Theorem 3 in \cite{KSzStudia} shows that our definition is consistent with the one given in the contractive case.

The following theorem illuminates the importance of quasianalytic operators, showing that the challenging \emph{Hyperinvariant Subspace Problem} can be reduced to this class in the asymptotically non-vanishing case.

\begin{thm}
\label{HSP for non-qa}
Let $T\in\L(\HH)$ be an absolutely continuous polynomially bounded operator, which is asymptotically non-vanishing: $\HH_0(T)\neq\HH$.
If $\, T$ is not quasianalytic, then it has a non-trivial hyperinvariant subspace.
\end{thm}

\noindent\textbf{Proof.}
By our assumption we can find non-zero vectors $x_1, x_2\in\HH$ so that $\o(T,x_2)$ is not contained in $\o(T,x_1)$, that is the set $\Delta= \o(T,x_2)\setminus\o(T,x_1)$ is of positive measure.
There exists a decreasing sequence $F= \{f_n\}_{n=1}^\i$ in $H^\i$ such that $x_1\in\HH_0(T,F)$ and $N_F=\Delta$; see the second part of the proof of Proposition \ref{bounded sequences}.
The latter relation implies by Theorem \ref{local stability} that $x_2\not\in\HH_0(T,F)$. 
Therefore, $\HH_0(T,F)$ is a proper hyperinvariant subspace of $T$.

\rightline{$\square$}

\medskip
The preceding proof immediately yields the following statement.

\begin{pro}
\label{two local residual sets}
Let $T_1\in\L(\HH_1)$ and $T_2\in\L(\HH_2)$ be a.c.\ polynomially bounded operators.
If $\o(T_2,x_2)$ is not contained in $\o(T_1,x_1)$, then there exists a decreasing sequence $F$ in $H^\i$ such that $x_1\in \HH_0(T_1,F)$ and $x_2\not\in\HH_0(T_2,F)$.
\end{pro}

\medskip

We proceed with some intertwining relations.

\begin{pro}
\label{intertwining relations}
Let $T_1\in\L(\HH_1)$ and $T_2\in\L(\HH_2)$ be a.c.\ polynomially bounded operators.
\begin{itemize}
\item[\textup{(a)}] If $\, Y\in\I(T_1,T_2)$, then $\o(T_1,x)\supset \o(T_2,Yx)$ holds for every $x\in\HH_1$.
\item[\textup{(b)}] If $\, \I(T_1,T_2)$ contains an injection $Y$, then $\pi(T_1)\supset\pi(T_2)$.
\item[\textup{(c)}] If both $\, \I(T_1,T_2)$ and $\, \I(T_2,T_1)$ contain injections, then $\pi(T_1)=\pi(T_2)$.
\item[\textup{(d)}] If $\, \I(T_1,T_2)$ contains a transformation $Z$ with dense range, then $\o(T_1)\supset\o(T_2)$.
\item[\textup{(e)}] If $\, T_1 \sim T_2$, then $\pi(T_1)=\pi(T_2)$ and $\o(T_1)=\o(T_2)$.
\item[\textup{(f)}] If $T_1 \sim T_2$ and $T_1$ is quasianalytic, then  $T_2$ is also quasianalytic.
\end{itemize}
\end{pro}

\noindent\textbf{Proof.}
The equations $Yf_n(T_1)= f_n(T_2)Y\; (n\in\N)$ show that $Y\HH_0(T_1,F)\subset \HH_0(T_2,F)$ holds, for any decreasing sequence $F= \{f_n\}_{n=1}^\i$ in $H^\i$.
In view of  Proposition \ref{two local residual sets} we infer that $\o(T_1,x)\supset\o(T_2,Yx)$ for every $x\in\HH_1$.
Assuming that $Y\in\I(T_1,T_2)$ is injective, this relation yields 
\[\pi(T_2)= \wedge\{\o(T_2,y) : 0\neq y\in\HH_2\} \subset \wedge\{\o(T_1,x) : 0\neq x\in\HH_1\}=\pi(T_1).\]
Let us assume now that $Z\in\I(T_1,T_2)$ has dense range.
Let $(X_j,U_j)$ be a unitary asymptote of $T_j\; (j=1,2)$.
Since $X_2Z\in\I(T_1,U_2)$, there exists a unique $\wt Z\in \I(U_1,U_2)$ such that $X_2Z= \wt Z X_1$.
It follows that
\[(\wt Z\K_1)\ba = (\wt Z \vee_{n\in\N}U_1^{-n}X_1\HH_1)\ba = \vee_{n\in\N} U_2^{-n} \wt Z X_1\HH_1= \vee_{n\in\N} U_2^{-n}X_2\HH_2=\K_2.\]
Considering the polar decomposition $\wt Z= W |\wt Z|$ we obtain that $W\in\I(U_1,U_2)$ is a coisometry, and so $U_2$ is unitarily equivalent to the restriction of $U_1$ to its reducing subspace $(\hbox{ker } W)^{\perp}$.
Therefore, the measurable support $\o(T_1)$ of the spectral measure of $U_1$ contains the measurable support $\o(T_2)$ of the spectral measure of $U_2$.

\rightline{$\square$}

\medskip
Quasianalicity determines the asymptotic behaviour of the operator.

\begin{pro}
\label{asymptotic behaviour}
If the a.c.\ polynomially bounded operator $T\in\L(\HH)$ is quasianalytic, then $T$ is of class $C_{10}: \; \HH_0(T)= \{0\}$ and $\HH_0(T^*)=\HH$.
\end{pro}

\noindent\textbf{Proof.}
Taking the decreasing sequence $F_0= \{\h^n\}_{n=1}^\i$ in $H^\i$, we infer that $\HH_0(T)= \HH_0(T,F_0)= \{0\}$ since $N_{F_0}\cap \pi(T)= \pi(T)\neq\emptyset$; see Proposition \ref{bounded sequences}.
Assuming that $\HH_0(T^*)\neq\HH$, let us consider a unitary asymptote $(X_*^*,U_*^*)$ of the adjoint $T^*$.
Then the a.c.\ unitary operator $U_*$ acts on a non-zero space $\K_*$, and $X_*\in\I(U_*,T)$ is a non-zero transformation.
The equations $X_*U_*= TX_*$ and $X_*=TX_*U_*^*$ imply that the subspace $\hbox{ ker}\,X_*$ is reducing for $U_*$, because of  $\hbox{ker}\,T\subset \HH_0(T)= \{0\}$.
Let us consider the restriction $U_{*,0}$ of $U_*$ to the non-zero reducing subspace $\K_* \ominus \hbox{ ker}\,X_*$.
Since $\I(U_{*,0},T)$ contains an injection, we infer by Proposition \ref{intertwining relations} that $\pi(U_{*,0}) \supset \pi(T)$.
However, it is clear that $\pi(U_{*,0})=\emptyset$ holds for the a.c.\ unitary operator $U_{*,0}$, and so $\pi(T)=\emptyset$, what is a contradiction.

\rightline{$\square$}

Now we examine the effect of taking orthogonal sums.

\begin{pro}
\label{orthogonal sums}
Let $T_j\in\L(\HH_j)$ be a.c.\ polynomially bounded operator acting on non-zero space, for $j=1,2$, and let us form $T= T_1 \op T_2 \in\L(\HH=\HH_1\op \HH_2)$.
\begin{itemize}
\item[\textup{(a)}] Then $T$ is also an a.c.\ polynomially bounded operator with \[\pi(T)= \pi(T_1)\cap \pi(T_2) \quad \hbox{ and } \quad \o(T)=\o(T_1)\cup\o(T_2).\]
\item[\textup{(b)}] $T$ is quasianalytic if and only if $T_1$ and $T_2$ are quasianalytic and $\pi(T_1)=\pi(T_2)$.
\end{itemize}
\end{pro}

\noindent\textbf{Proof.}
Setting $x=x_1\op x_2, \; y=y_1\op y_2\in\HH$, and $\m_1\in M(T_1,x_1,y_1), \; \m_2\in M(T_2,x_2,y_2)$, we have
\[\la u(T)x,y\ra = \la u(T_1)x_1,y_1\ra + \la u(T_2)x_2,y_2\ra = \int_\T u\, d(\m_1+\m_2)\]
for every $u\in \A(\T)$; hence $\m_1+\m_2\in M(T,x,y)$.
Therefore $T$ is an a.c.\  polynomially bounded operator.

Let $(X_j,U_j)$ be a unitary asymptote of $T_j$, for $j=1,2$.
We know that $(X,U)$ will be a unitary asymptote of $T$, where $X=X_1\op X_2$ and $U=U_1\op U_2$; see Corollary \ref{direct sum}.
Furthermore, if $E$ is the spectral measure of $U$, then $E|\K_j$ will be the spectral measure of $U_j\; (j=1,2)$.
Thus, for any $x=x_1\op x_2\in\HH$ we have $E_{Xx,Xx}= E_{X_1x_1,X_1x_1} + E_{X_2x_2,X_2x_2}$, whence $\o(T,x)= \o(T_1,x_1)\cup\o(T_2,x_2)$ can be derived.
From here the statement follows.

\rightline{$\square$}

\medskip
We conclude with the properties of invariant subspaces.

\begin{pro}
\label{restrictions}
Let $T\in\L(\HH)$ be an a.c.\ polynomially bounded operator, and let us assume that $\M\in \textup{ Lat}\,T$ is non-zero.
\begin{itemize}
\item[\textup{(a)}] Then $T|\M$ is also an a.c\ polynomially bounded operator with 
\[\pi(T)\subset \pi(T|\M) \subset \o(T|\M)\subset \o(T).\]
\item[\textup{(b)}] If $\, T$ is quasianalytic, then $T|\M$ is also quasianalytic and $\pi(T|\M)= \pi(T)$.
\end{itemize}
\end{pro}

\noindent\textbf{Proof.}
It is clear that $M(T|\M,x,y)= M(T,x,y)$ holds for every $x,y\in\M$.
Hence, $T|\M$ is an a.c.\ polynomially bounded operator.
We know that if $(X,U)$ is a unitary asymptote of $T$, then $(X|\M,U|\wt\M)$ is a unitary asymptote of $T|\M$; see Corollary \ref{restriction}.
Thus, $\o(T|\M,x)= \o(T,x)$ holds, for every $x\in\M$.

\rightline{$\square$}

\medskip
We recall that if $T\in\L(\HH)$ is an a.c.\  contraction with $\o(T)=\T$, then there exists $\M\in \textup{ Lat}\,T$ such that $T|\M$ is similar to the simple unilateral shift $S\in\L(H^2),\; Sf=\h f$; even more, these shift-type invariant subspaces span the whole space $\HH$ (see Theorem IX.3.6 in \cite{NFBK} and \cite{KerchyJFA}).

\begin{que}
\label{full residual set}
Let $T\in\L(\HH)$ be an a.c.\  polynomially bounded operator with $\o(T)=\T$ (or $\pi(T)=\T$). Does there exist a non-zero $\M\in \textup{ Lat}\,T$ such that $T|\M$ is similar to a contraction?
\end{que}

\medskip
If $\,\pi(T)=\T$, then such an $\M$ clearly contains a subspace $\M'\in \textup{ Lat}\,T$, where $T|\M'$ is similar to $S$.
It would be interesting to know whether Pisier's construction provides completely non-contractive (i.e., having no $\M\in \textup{ Lat}\,T \setminus\{0\}$ such that $T|\M$ is similar to a contraction) absolutely continuous polynomially bounded operators in the case, when $\o(T)=\T$ or $\pi(T)=\T$.

\medskip
\noindent
\textsc{L. K\'erchy}, Bolyai Institute, University of Szeged, Aradi v\'ertan\'uk tere 1, H-6720 Szeged, Hungary; \emph{e-mail}: kerchy@math.u-szeged.hu

\end{document}